\documentclass[a4paper,11pt]{article}
\usepackage[utf8]{inputenc}
\usepackage[T1]{fontenc}
\usepackage[bitstream-charter]{mathdesign}

\usepackage{cite}
\usepackage{amsmath}
\usepackage{amssymb}
\usepackage{amsthm}
\usepackage{xcolor}
\usepackage{graphicx}
\usepackage{url}
\usepackage{enumitem}
\usepackage{accents}
\usepackage{marvosym}
\usepackage{authblk}
\usepackage{vmargin}
\setmarginsrb{2cm}{2cm}{2cm}{2cm}{0mm}{0mm}{0mm}{10mm}
\theoremstyle{definition}

\newtheorem{theorem}{Theorem}

\newtheorem{corollary}{Corollary}
\newtheorem{lemma}{Lemma}
\newtheorem{remark}{Remark}
\usepackage{hyperref}
\hypersetup{colorlinks=true,allcolors=blue}
\DeclareFontShape{OMX}{cmex}{m}{n}{%
	<-7.5>cmex7
	<7.5-8.5>cmex8
	<8.5-9.5>cmex9
	<9.5->cmex10
}{}
\DeclareSymbolFont{largesymbols}{OMX}{cmex}{m}{n}

\title{\vspace*{-1.5cm} \bfseries Unraveling the complexity of inverting the Sturm-Liouville boundary value problem to its canonical form} 
\author[1]{N. Karjanto\thanks{\Letter: \url{natanael@skku.edu}\href{https://orcid.org/0000-0002-6859-447X}{\includegraphics[scale=0.08]{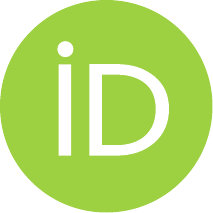}}}}
\author[2]{P. Sadhani\,}
\affil[1]{Department of Mathematics, University College, Natural Science Campus, Sungkyunkwan University\protect \\ 2066~Seobu-ro, Jangan-gu, Suwon~16419, Gyeonggi-do, Republic of Korea}
\affil[2]{Department of Computer Science, University of Oxford, Wellington Square, Oxford OX1 2JD, United Kingdom}

\date{\vspace*{-0.5cm} \footnotesize Last updated \today}

\begin{document}
\maketitle

\begin{abstract}
\noindent	
The Sturm-Liouville boundary value problem (SLBVP) stands as a fundamental cornerstone in the realm of mathematical analysis and physical modeling. Also known as the Sturm-Liouville problem (SLP), this paper explores the intricacies of this classical problem, particularly the relationship between its canonical and Liouville normal (Schrödinger) forms. While the conversion from canonical to Schrödinger form using Liouville's transformation is well-known in the literature, the inverse transformation seems to be neglected. Our study attempts to fill this gap by investigating the inverse of Liouville's transformation, that is, given any SLP in the Schrödinger form with some invariant function, we seek the SLP in its canonical form. By examining closely the second Paine-de Hoog-Anderson (PdHA) problem, we argue that retrieving the SLP to its canonical form can be notoriously difficult and even impossible to achieve in its exact form. Finding the inverse relationship between the two independent variables seems to be the main obstacle. We confirm this claim by considering four different scenarios depending on the potential and density functions that appear in the corresponding invariant function. In the second PdHA problem, this invariant function takes a reciprocal quadratic binomial form. In some cases, the inverse Liouville's transformation produces an exact expression for the SLP in its canonical form. In other situations, however, while an exact canonical form is not possible to obtain, we have successfully derived the SLP in its canonical form asymptotically. \\

\end{abstract}

\section{Introduction}

The Sturm-Liouville boundary value problem (SLBVP) stands as a venerable and foundational topic in the realm of mathematical analysis and scientific inquiry. The SLBVP is often referred to SLP (Sturm-Liouville problem) in the literature, and we will follow this convention in this article. The SLP is a specific type of second-order linear ordinary differential equation (ODE) problem that arises in various areas of mathematics and physics, particularly in the study of partial differential equations (PDEs) and the eigenvalue problems associated with them. It is named after the mathematicians who pioneered the study, that is, Jacques Charles François Sturm (1803--1855) and Joseph Liouville (1809--1882). 

The origin of SLP traces its history back to the first half of the 19th century, when both Sturm and Liouville published a sequence of papers on second-order linear ODEs that include BVPs between 1836 and 1837. The study of ODEs before their period was predominantly limited to searching for solutions in terms of analytic expressions. Sturm and Liouville were among the first mathematicians who recognized the limitations of such approach and discerned the advantage in investigating properties of solutions directly from the ODEs even in the absence of closed forms for solutions~\cite{zettl2010sturm,luetzen1984sturm,luetzen2012joseph}. 

Thanks to the advanced progress in computational tools, it only emerged during the second half of the 20th century that the Sturm-Liouville theory and its applications attracted significant attention. Even though the problem has become a classic, there are still dozens of research papers on the topic published annually. The mathematical framework of SLP remains of paramount importance across various disciplines of mathematics, physics, and engineering. The SLP serves as a versatile and elegant tool for understanding and analyzing the behavior of second-order linear ODEs, with diverse applications in quantum mechanics, heat transfer, structural dynamics, vibrations of mechanical systems, and many more. Solving these problems provides the eigenvalues and eigenfunctions associated with these systems, with important physical and mathematical interpretations~\cite{pryce1993numerical,bailey2001sleighn2,amrein2005sturm,gwaiz2008sturm,agarwal2009ordinary,teschl2012ordinary,haberman2013applied,guenther2019sturm,kravchenko2020direct,masjed-jamei2020special,zettl2021recent}.

Unless otherwise mentioned, we will focus on the classical SLP in this article, which refers to the original and foundational formulation of the problem with specific properties and assumptions. The term \emph{classical} distinguishes this foundational formulation from more generalized or extended versions of the SLP. The classical SLP serves as a fundamental model, and its properties and solutions have also been extensively studied and applied in various areas of mathematics and physics. Some key characteristics of the classical SLP include second-order linear ODEs, self-adjoint (Hermitian) differential operators, real-valued eigenvalues, and specific (e.g., Dirichlet, Neumann, Robin) boundary conditions. The generalizations of SLP may encompass fractional differential equations~\cite{almdallal2009efficient,klimet2013fractional,zayernouri2013fractional}, non-self-adjoint operators~\cite{veliev2007non,albeverio2008on,chanane2007computing}, complex-valued eigenvalues~\cite{behrndt2008accumulation,xie2013non,levitin2015accumulation}, or different (e.g., time-dependent, eigenparameter-dependent, nonlocal, etc.) boundary conditions~\cite{sager1984the,binding1994sturm,stikonas2007the}. 

Our study investigates the transformation between the Sturm-Liouville problem (SLP) in its canonical form and the Liouville normal (Schrödinger) form. This conversion, known as Liouville's transformation/reduction, involves transforming both the dependent and independent variables~\cite{pryce1993numerical,adam2017rays,birkhoff1989ordinary,everitt1982on}. While the transformation itself is established, the inverse Liouville's transformation, which retrieves the SLP from its Schrödinger form, remains largely unexplored. This lack of attention might be explained by the assumption that reversing the steps of the transformation would be trivial. However, as we will see in the case of the generalized second Paine-de Hoog-Anderson (PdHA) problem, inverting the SLP to its canonical form can be notoriously difficult even after solving the associated ODE for the invariant function.

This work fills a gap in the understanding of the SLP transformation between its two forms. We elucidate the process of conversion by considering the second PdHA problem, also known in the literature as the second Paine problem. We focus on a generalized version of the invariant function, focusing on the reciprocal quadratic power of the binomial term as originally presented in their 1981 paper~\cite{paine1981correction}. Interestingly, depending on the combinations of potential and density functions, the inversion process can range from relatively straightforward to extremely difficult. As a consequence, the SLP can be retrieved to its canonical form in either case: with an exact expression of the so-called $p$-function for the former case, or as an asymptotic approximation for the latter, where we propose a novel technique using asymptotic expansion. Therefore, we are able to obtain the SLP in its canonical form in either case, either with an exact expression or as an asymptotic approximation. This constitutes our primary contribution.

In her PhD thesis, Ledoux discussed a new class of numerical methods in detail for solving the SLP in both canonical and Schrödinger forms. Using the constant and line perturbation methods, which fall under the category of piecewise perturbation methods, these techniques demonstrated efficiency, accuracy, and stability~\cite{ledoux2007study}. She briefly touched on Liouville’s transformation, and in the appendix she provided some functions associated with the second Paine problem in the canonical form, referring to~\cite{ixaru1999slcpm}. However, no derivation was provided on how to obtain these functions. Therefore, we will delve deeper into this issue and provide the missing derivation at the end of Section~\ref{SectionPaine-n2}. Interestingly, despite Everitt providing a catalogue of nearly 60~examples of SLP in either canonical or Liouville normal forms, both the first and second Paine problems were conspicuously absent from the list~\cite{everitt2005a}. 

Transforming the SLP from its Liouville normal form to the canonical form has been studied in the context of perturbed potential temperature fields in atmospheric boundary layers. While close convergence was observed between asymptotic solutions using the WKB method and numerical simulations, the study only focused on a simplified case with constant density and zero potential~\cite{karjanto2022perturbed}. Another study investigated the lowest-order eigenvalue under both Dirichlet and Neumann boundary conditions to further understand the second generalized Paine problem~\cite{karjanto2022modified}. Using a method for estimating the lowest eigenvalue that incorporates the localized landscape and potential functions~\cite{arnold2019computing,filoche2016universal}, the study found that the estimates tend to overshoot the actual values but exhibit excellent qualitative agreement for the second Paine problem with Dirichlet boundary conditions. However, discrepancies were observed for Neumann boundary conditions with specific combinations of small constant in the binomial term and higher values of the denominator.

This article is organized as follows. Section~\ref{SectionSLBVP} discusses the transformation of the SLP from its canonical form to the Schr\"odinger form and its inverse. The so-called Liouville's transformation involves both the independent and dependent variables of the {SLP}. Section~\ref{SectionPaine-n2} focuses on the second PdHA problem, which generalizes the corresponding invariant function to include arbitrary positive constants while maintaining its reciprocal quadratic power. We consider four different scenarios, each discussed in its own subsection. The first case examines a combination of a nonzero constant density function and a vanishing potential function. The second case explains the combination of a nonzero constant potential and quadratic density functions. The third case delves into the situation where both the potential and density functions are nonzero constants. The final case discuses when the transformation function depends reciprocally linear on the Schr\"odinger variable. Finally, Section~\ref{SectionConclude} concludes our discussion and provide further remarks and future directions.

\section{Sturm-Liouville boundary value problem}	\label{SectionSLBVP}

In this section, we reconfigure the SLP in the canonical form to the one in the Schr\"odinger form using a transformation that involve both the independent and dependent variables, the so-called Liouville's transformation. We also verify that, using the inverse transformation, the SLP in the canonical form can be recovered from the SLP in the Liouville normal form.

Consider the general form of the classical Sturm-Liouville eigenvalue problem, written in the {\color{blue} \emph{canonical form}} with eigenvalue $\lambda$ and the corresponding eigenfunction $u\left(x\right)$:
\begin{equation}
-\frac{d}{dx}\left(p \left(x \right) \frac{du}{dx} \right) + q \left(x \right) u = \lambda \, r \left(x \right) u, \qquad a < x < b, \quad \text{or} \quad x \in \mathbb{R}.		\label{EquSLP}
\end{equation}
The regular boundary conditions are imposed at the endpoints:
\begin{equation*}
\delta_0 \, u\left(a \right) - \delta_1 p\left(a \right) \frac{du}{dx}\left(a \right) = 0, \qquad \qquad \gamma_0 \, u\left(b \right) - \gamma_1 p \left(b \right) \frac{du}{dx}\left(b \right) = 0,
\end{equation*}
where $\delta_0$ and $\delta_1$ are real and not both zero, and similarly, $\gamma_0, \gamma_1 \in \mathbb{R}$, $\gamma_0^2 + \gamma_1^2 > 0$. The SLP~\eqref{EquSLP} is regular if both $a$ and $b$ are finite and the functions $p$, $q$, and $r$ are piecewise continuous on $[a,b]$, usually not zero, in particular, $p, r > 0$, but they may take constant values. The function $r = r\left(x\right)$ is called the \emph{weight} or \emph{density} function, whereas $q = q\left(x\right)$ is often referred to as the \emph{potential} function by physicists and chemists. 

The following lemma considers a special case when both $p$ and $r$ are nonzero constants.
\begin{lemma}
When both $p$ and $r$ are constants but nonzero, the ODE~\eqref{EquSLP} can be expressed in the Liouville normal form using the change of variable $t = \eta x$, where $\eta = \sqrt{r/p} \neq 0$:
\begin{equation}
-\frac{d^2u}{dt^2} + Q\left(t\right) \, u = \lambda u, \qquad \alpha < t < \beta, \quad \text{where} \quad Q\left(t\right) = \frac{pq}{w^2}.		\label{EquLemma0}
\end{equation}
\end{lemma}

\begin{proof}
Because $p$ is constant and $p \neq 0$, we can pull it out from the derivative sign and divide ODE~\eqref{EquSLP} by $p$:
\begin{equation}
-\frac{d^2u}{dx^2} + \frac{q}{p} \, u = \lambda \frac{r}{p} \, u.		\label{Equconstantpw}
\end{equation}
Using the transformation $t = \eta x$, we have
\begin{equation*}
\frac{du}{dx} = \eta \frac{du}{dt}, \qquad \text{and} \qquad
\frac{d^2u}{dx^2} = \eta^2 \frac{d^2u}{dt^2}.
\end{equation*}
The ODE~\eqref{Equconstantpw} becomes
\begin{equation*}
-\eta^2 \frac{d^2u}{dt^2} + \frac{q}{p} \, u = \lambda \, \eta^2 \, u.
\end{equation*}
By dividing with $\eta^2 \neq 0$, we obtain the desired form~\eqref{EquLemma0}, where
\begin{equation*}
Q\left(t\right) = \frac{q}{\eta^2 \, p} = \frac{pq}{r^2}.
\end{equation*}
All expressions for the functions $p$, $q$, and $r$ are to be understood as $f\left(x \right) = f \left(t/\eta \right)$, where $f = p, q$, or $r$. This completes the proof.
\end{proof}

\begin{remark}
In the context of quantum mechanics, the ODE~\eqref{EquLemma0} represents the time-dependent Schr\"odinger equation, where $Q$ denotes the potential function and the eigenvalue $\lambda$ corresponds to the energy level. However, in this article, we will deviate from common physics terminology and instead use the term ``potential'' to the function $q$ in~\eqref{EquSLP} but not $Q$.
\end{remark}

The following lemma shows what happens to SLP~\eqref{EquSLP} when we transform its independent variable.
\begin{lemma}   \label{lemmaindependent}
The eigenvalue problem~\eqref{EquSLP} can be converted to another ODE in the following form by transforming the independent variable $x$ to $t$, where $x = x\left(t\right)$:
\begin{equation}
-\frac{d}{dt}\left(\frac{p}{\dot{x}} \, \frac{du}{dt} \right) + q \dot{x} \, u = \lambda \, r \dot{x} \, u,		\label{EquLemma1}
\end{equation}
where the dot represents the derivative with respect to $t$, that is, $\dot{x} = dx/dt$, and is assumed to take one sign on the open interval $\alpha < t < \beta$, that corresponds to $a < x < b$ in the original variable through a transformation. Furthermore, $u = u\left(x \left(t \right) \right)$, and $u/\dot{x} = u \left(x \left(t \right)\right)/\dot{x}\left(t\right)$. Other functions also follow a similar convention.
\end{lemma}

\begin{proof}
It can be easily worked out using the fact that the differential operator $d/dx = \left(dt/dx \right) \, d/dt = \left(1/\dot{x} \right) \, d/dt$, and thus $du/dx = \left(1/\dot{x} \right) \, du/dt$. Substituting these expressions to ODE~\eqref{EquSLP}, we obtain
\begin{equation}
-\frac{1}{\dot{x}} \frac{d}{dt} \left( \frac{p}{\dot{x}} \frac{du}{dt} \right) + q \, u = \lambda \, r \, u.		\label{EquLem1indevar}
\end{equation}
Multiplying both sides of~\eqref{EquLem1indevar} with $\dot{x}$ yields the desired expression~\eqref{EquLemma1}:
\begin{equation*}
-\frac{d}{dt}\left(\frac{p}{\dot{x}} \frac{du}{dt} \right) + q \dot{x} \, u = \lambda \, r \dot{x} \, u.
\end{equation*}
The proof is completed.
\end{proof}

The following lemma demonstrates how eigenvalue problem~\eqref{EquSLP} transforms to another ODE when we introduce a new dependent variable.
\begin{lemma}   \label{lemmadependent}
The eigenvalue problem~\eqref{EquSLP} can be converted to another ODE in the following form by transforming the dependent variable $u$ of the form $u \left(x \right) = w \left(x \right) \, v\left(x \right)$, with $w \left(x \right)$ is a given function:
\begin{equation}
-\frac{d}{dx}\left(p w^2 \, \frac{dv}{dx}\right) + \left[q w^2 - w \frac{d}{dx} \left(p \frac{dw}{dx} \right) \right] v = \lambda \, r w^2 \, v.		\label{EquLemma2}
\end{equation}
\end{lemma}

\begin{proof}
By applying the product rule to $u$, we have
\begin{align*}
\frac{du}{dx} &= \frac{dw}{dx} v + w \frac{dv}{dx}, \\
p\frac{du}{dx} &= pv \frac{dw}{dx} + p w \frac{dv}{dx}, \\
-\frac{d}{dx} \left(p \frac{du}{dx} \right) &= -\frac{d}{dx} \left( pv \frac{dw}{dx} \right) - \frac{d}{dx} \left( p w \frac{dv}{dx} \right).
\end{align*}
Substituting these expressions to ODE~\eqref{EquSLP} yields
\begin{equation}
-\frac{d}{dx} \left( p w \frac{dv}{dx} \right) + q w \, v - \frac{d}{dx} \left( pv \frac{dw}{dx} \right) = \lambda \, r w \, v.		\label{EquLem2a}
\end{equation}
Multiplying~\eqref{EquLem2a} with the function $w$ results
\begin{equation}
- w \frac{d}{dx} \left( p w \frac{dv}{dx} \right) + q w^2 \, v - w \frac{d}{dx} \left( pv \frac{dm}{dx} \right) = \lambda \, r w^2 \, v.	\label{EquLem2b}
\end{equation}
Expanding the third term on the left-hand side of~\eqref{EquLem2b}, we obtain
\begin{equation*}
- w \frac{d}{dx} \left( p w \frac{dv}{dx} \right) + q w^2 \, v - w v\frac{d}{dx} \left( p \frac{dw}{dx} \right) {\color{blue} \, - \, wp \frac{dv}{dx} \frac{dw}{dx}} = \lambda \, r w^2 \, v.
\end{equation*}
Because
\begin{equation*}
\frac{d}{dx} \left(p w^2 \frac{dv}{dx} \right) = w \frac{d}{dx} \left(pw \frac{dv}{dx} \right) {\color{red} \, + \, pw \frac{dw}{dx} \frac{dv}{dx}},
\end{equation*}
we observe that two terms ({\color{blue} blue} and {\color{red} red}) cancel each other:
\begin{equation}
-\frac{d}{dx} \left( p w^2 \frac{dv}{dx} \right) {\color{red} + \, pw \frac{dw}{dx} \frac{dv}{dx}} + q w^2 \, v - w v \frac{d}{dx} \left( p \frac{dw}{dx} \right) {\color{blue} - wp \frac{dv}{dx} \frac{dw}{dx} } = \lambda \, r w^2 \, v.	 \label{EquLem2c}
\end{equation}
Rearranging the remaining terms of~\eqref{EquLem2c}, we obtain the desired ODE~\eqref{EquLemma2}:
\begin{equation*}
-\frac{d}{dx}\left(p w^2 \, \frac{dv}{dx}\right) + \left[q w^2 - w \frac{d}{dx} \left(p \frac{dw}{dx} \right) \right] v = \lambda \, rw^2 \, v.
\end{equation*}
This completes the proof.
\end{proof}

Lemma~\ref{lemmacombined} combines both Lemma~\ref{lemmaindependent} and Lemma~\ref{lemmadependent}, and reveals what kind of eigenvalue problem will be obtained when both independent and dependent variables are transformed simultaneously.
\begin{lemma}		\label{lemmacombined}
By combining both the transformations for independent and dependent variables, that is, $x = x\left(t\right)$ and $u\left(x\right) = w\left(x\right) \, v\left(x\right)$, we obtain another transformed Sturm-Liouville ODE:
\begin{equation}
-\frac{d}{dt} \left(P \left(t\right) \frac{dv}{dt} \right) + Q\left(t\right) \, v = \lambda \, R\left(t\right) v,		\label{EquLemma3}
\end{equation}
where
\begin{align*}
P\left(t\right) &= \frac{pw^2}{\dot{x}}, \\
Q\left(t\right) &= \left[qw - \frac{1}{\dot{x}} \frac{d}{dt} \left(\frac{p}{\dot{x}} \frac{dw}{dt} \right) \right] \, w \dot{x}, \\
R\left(t\right) &= r w^2 \dot{x}.
\end{align*}
\end{lemma}

\begin{proof}
Using the results from Lemma~\ref{lemmaindependent} and Lemma~\ref{lemmadependent}, we can write the ODE~\eqref{EquSLP} as follows:
\begin{equation}
-\frac{1}{\dot{x}}\frac{d}{dt} \left(\frac{p w^2}{\dot{x}} \, \frac{dv}{dt}\right) + \left[q w - \frac{1}{\dot{x}} \frac{d}{dt} \left(\frac{p}{\dot{x}} \frac{dw}{dt} \right) \right] w \, v = \lambda \, r w^2 \, v. 	\label{EquLem3a}
\end{equation}
Multiplying~\eqref{EquLem3a} with $\dot{x}$, we obtain the desired result~\eqref{EquLemma3}. This completes the proof.
\end{proof}

The following theorem describes how Liouville's transformation convert the SLP~\eqref{EquSLP} in the canonical form to another SLP in the Liouville normal form.
\begin{theorem}[Liouville's transformation~\cite{liouville1837second}]
The Sturm-Liouville problem in the {\color{blue} \emph{canonical form}} with eigenvalue~$\lambda$ and the corresponding eigenfunction~$u \left(x\right)$, with regular boundary conditions:
\begin{align*}
-\frac{d}{dx} \left[p \left(x\right) \frac{du}{dx} \right] + q\left(x\right) \, u = \lambda \, r\left(x\right) \, u, & \qquad \qquad a < x < b, \\
\delta_0 \, u\left(a \right) - \delta_1 \, p\left(a \right) \frac{du}{dx}\left(a\right) = 0, & \qquad \qquad \gamma_0 \, u\left(b\right) - \gamma_1 \, p \left(b\right) \frac{du}{dx}\left(b\right) = 0,
\end{align*}
where $\delta_0$ and $\delta_1$ are real and not both zero, $\gamma_0$ and $\gamma_1$ are also similarly conditioned, can be converted into the {\color{blue} \emph{Liouville normal (Schr\"odinger) form}} by performing {\color{blue} \emph{Liouville's transformation}}
\begin{align*}
-\frac{d^2v}{dt^2} + I\left(t\right) \, v = \lambda \, v,& \qquad \qquad \alpha < t < \beta,\\
\delta_2 \, v\left(\alpha\right) - \delta_1 \, P\left(\alpha\right) \frac{dv}{dt}\left(\alpha\right) = 0, & \qquad \qquad
\gamma_2 \, v\left(\beta\right)  - \gamma_1 \, P\left(\beta\right)  \frac{dv}{dt}\left(\beta\right)  = 0,
\end{align*}
where
\begin{align*}
\delta_2 = \left. \left(\delta_0 w^2 - \delta_1 \, p w \frac{dw}{dx} \right) \right|_{x = a} \qquad \text{and} \qquad
\gamma_2 = \left. \left(\gamma_0 w^2 - \gamma_1 \, p w \frac{dw}{dx} \right) \right|_{x = b}.
\end{align*}
Here, $I$ is the corresponding {\color{blue} \emph{invariant function}} of the SLP~\eqref{EquSLP}:
\begin{equation}
I\left(t\right) = \frac{q}{r} + w \frac{d^2}{dt^2} \left(\frac{1}{w}\right),		\label{Equ-invariant-function}
\end{equation}
and the {\color{blue} \emph{Liouville transformation}} is given by
\begin{equation}
t = \int \sqrt{\frac{r}{p}} \, dx, \qquad w = \left(p r \right)^{-1/4}, \qquad \text{and} \qquad u\left(x\right) = w\left(x\right) \, v\left(x\right).		\label{Equ-Liouville-transformation}
\end{equation}
\end{theorem}

\begin{proof}
Using the independent variable transformation $\left(x \leftrightarrow t \right)$, we have
\begin{equation*}
\frac{dt}{dx} = \frac{1}{dx/dt} = \frac{1}{\dot{x}} = \sqrt{\frac{r}{p}}, \qquad \qquad \text{or} \qquad \qquad \dot{x} = \sqrt{\frac{p}{r}}.
\end{equation*}
Because $w^2 = 1/\sqrt{pr}$ or $1/w^2 = \sqrt{pr}$, we also have
\begin{equation*}
\frac{p}{\dot{x}} \frac{dw}{dt} = p \sqrt{\frac{r}{p}} \frac{dw}{dt} = \sqrt{rp} \frac{dw}{dt} = \frac{1}{w^2} \frac{dw}{dt}.
\end{equation*}
Using Lemma~\ref{lemmacombined}, we observe that
\begin{align*}
P\left(t\right) &= \frac{pw^2}{\dot{x}} = \frac{p}{\sqrt{pr}} \sqrt{\frac{r}{p}} = 1, \\
Q\left(t\right) &= qw^2 \dot{x} - w \frac{d}{dt} \left(\frac{p}{\dot{x}} \frac{dw}{dt} \right)
      = \frac{q}{\sqrt{pr}} \sqrt{\frac{p}{r}} + w \frac{d}{dt} \left(-\frac{1}{w^2} \frac{dw}{dt} \right) \\
     &= \frac{q}{r} + w \frac{d}{dt} \left[ \frac{d}{dt} \left(\frac{1}{w}\right) \right] 
      = \frac{q}{r} + w \frac{d^2}{dt^2} \left(\frac{1}{w}\right) = I\left(t\right) \\
R\left(t\right) &= r w^2 \dot{x} = \frac{r}{\sqrt{pr}} \sqrt{\frac{p}{r}} = 1.      
\end{align*}
Hence, the ODE~\eqref{EquLemma3} becomes
\begin{equation*}
-\frac{d^2v}{dt^2} + I(t) \, v = \lambda \, v, \qquad \qquad \alpha < t < \beta.
\end{equation*}
For the boundary conditions, multiply both of them with $w$ and use the fact that
\begin{equation*}
\frac{du}{dx} = \frac{dw}{dx} v + w \frac{dv}{dx} = \frac{dw}{dx} v + \frac{w}{\dot{x}} \frac{dv}{dt}.
\end{equation*}
By substituting the relevant boundaries, that is $x = a$ or $t = \alpha$, we obtain
\begin{align*}
\delta_0 \, w^2 v - \delta_1 \, pw \left(\frac{dw}{dx} v + \frac{w}{\dot{x}} \frac{dv}{dt} \right) &= 0, \\
\left(\delta_0 \, w^2 - \delta_1 \, pw \frac{dw}{dx}\right) v - \delta_1 \, \frac{pw^2}{\dot{x}} \frac{dv}{dt} &=0.
\end{align*}
By taking
\begin{equation*}
\delta_2 = \left. \left(\delta_0 w^2 - \delta_1 \, p w \frac{dw}{dx} \right) \right|_{x = a} \qquad \text{and} \qquad
P\left(\alpha\right) = \frac{pw^2}{\dot{x}}\Bigg|_{t = \alpha}
\end{equation*}
we obtain the boundary condition at $t = \alpha$. A similar argument can be reached for the boundary conditions at $t = \beta$. This completes the proof.
\end{proof}

The following corollary discuss the opposite process, that is, given the SLP in the Liouville normal form, we can invert it back to the SLP in the canonical form.
\begin{corollary}[Inverse Liouville's transformation]
The SLP expressed in its Sch\"odinger form	
\begin{align}
-\frac{d^2v}{dt^2} + I\left(t\right) \, v = \lambda \, v,& \qquad \qquad \alpha < t < \beta,	\label{EquCorol1A}	\\
\delta_2 \, v\left(\alpha\right) - \delta_1 \, P\left(\alpha\right) \frac{dv}{dt}\left(\alpha\right) = 0, & \qquad \qquad
\gamma_2 \, v\left(\beta\right)  - \gamma_1 \, P\left(\beta\right)  \frac{dv}{dt}\left(\beta\right)  = 0,
\end{align}
can be inverted back to the same eigenvalue problem in its canonical form
\begin{align}
-\frac{d}{dx} \left[p\left(x\right) \frac{du}{dx} \right] + q\left(x\right) \, u = \lambda r\left(x\right) \, u, & \qquad \qquad a < x < b, 	\label{EquCorol1B}\\		
\delta_0 \, u\left(a\right) - \delta_1 \, p\left(a\right) \frac{du}{dx}\left(a\right) = 0, & \qquad \qquad \gamma_0 \, u\left(b\right) - \gamma_1 \, p\left(b\right) \frac{du}{dx}\left(b\right) = 0,
\end{align}
using the identical Liouville's transformation.
\end{corollary}

\begin{proof}
We use the following relationships for the first-order and second-order derivative operators, respectively:
\begin{equation*}
\frac{d}{dt} = \dot{x} \frac{d}{dx} = \frac{dx}{dt} \, \frac{d}{dx} \qquad \qquad \text{and} \qquad \qquad
\frac{d^2}{dt^2} = \ddot{x} \frac{d}{dx} + \dot{x}^2 \frac{d^2}{dx^2} = \frac{d^2x}{dt^2} \, \frac{d}{dx} + \left(\frac{dx}{dt}\right)^2 \, \frac{d^2}{dx^2}.
\end{equation*}
Substitute these expressions to ODE~\eqref{EquCorol1A} yields
\begin{equation}
-\frac{d^2x}{dt^2} \, \frac{dv}{dx} - \left(\frac{dx}{dt}\right)^2 \, \frac{d^2v}{dx^2} + I(x) \, v = \lambda \, v.		\label{EquCorSub}
\end{equation}
Because $v(x) = u(x)/w(x)$, the first and second derivatives of $v$ can be expressed as follows, respectively:
\begin{align*}
\frac{d}{dx} \left(\frac{u}{w}\right) &= \frac{1}{w} \frac{du}{dx} - \frac{u}{w^2} \frac{dw}{dx} \\
\frac{d^2}{dx^2} \left(\frac{u}{w}\right) &= \frac{d}{dx} \left(\frac{1}{w} \frac{du}{dx} \right) - \frac{d}{dx} \left(\frac{u}{w^2} \frac{dw}{dx} \right) \\
&= \frac{1}{w} \frac{d^2u}{dx^2} - \frac{1}{w^2} \frac{dw}{dx} \frac{du}{dx} - \frac{1}{w^2} \frac{dw}{dx} \frac{du}{dx} + \frac{2}{w^3} \left(\frac{dw}{dx}\right)^2 u - \frac{1}{w^2} \frac{d^2w}{dx^2} u \\
&= \frac{1}{w} \frac{d^2u}{dx^2} - \frac{2}{w^2} \frac{dw}{dx} \frac{du}{dx} + \frac{2}{w^3} \left(\frac{dw}{dx}\right)^2 u - \frac{1}{w^2} \frac{d^2w}{dx^2} u.
\end{align*}
Substitute these expressions to ODE~\eqref{EquCorSub} and multiply it with $w$, we obtain
\begin{equation}
\begin{aligned}
-\left(\frac{dx}{dt}\right)^2 \, \frac{d^2u}{dx^2} + \left[\frac{2}{w} \frac{dw}{dx} \left(\frac{dx}{dt}\right)^2 - \frac{d^2x}{dt^2} \right]\, \frac{du}{dx} & \\ + 
\left\{I\left(x\right) - \left[\frac{2}{w^2} \left(\frac{dw}{dx}\right)^2 - \frac{1}{w} \frac{d^2w}{dx^2} \right] \, \left(\frac{dx}{dt}\right)^2 + \frac{1}{w} \frac{dw}{dx} \frac{d^2x}{dt^2} \right\} \, u &= \lambda \, u. 		\label{EquCorSubLong}
\end{aligned}
\end{equation}
Consider again ODE~\eqref{EquCorol1B} from the SLP in the canonical form, where we have now divided it by $r$:
\begin{equation}
-\frac{p}{r} \, \frac{d^2u}{dx^2} - \frac{1}{r} \frac{dp}{dx} \, \frac{du}{dx} + \frac{q}{r} \, u = \lambda u.		\label{EquCorSubShort}
\end{equation}
We compare these two ODEs, those are, Equations~\eqref{EquCorSubLong} and~\eqref{EquCorSubShort}. Because they should be identical, we obtain the following relationship:
\begin{equation*}
\left(\frac{dx}{dt}\right)^2 = \frac{p}{r} \qquad \text{or} \qquad \frac{dx}{dt} =  \pm \sqrt{\frac{p}{r}},
\end{equation*}
Equating the coefficient for $du/dx$, we observe that
\begin{align*}
\frac{2}{w} \frac{dw}{dx} r \left(\frac{dx}{dt}\right)^2 - r \frac{d^2x}{dt^2} &= - \frac{dp}{dx}, \\
\frac{2}{w} \frac{dw}{dx} p - \frac{1}{2} r \frac{d}{dx} \left(\frac{p}{r} \right) &= - \frac{dp}{dx}, \\
\frac{2p}{w} \frac{dw}{dx} + \frac{dp}{dx} - \frac{1}{2} r \left(\frac{1}{r} \frac{dp}{dx} - \frac{p}{r^2} \frac{dr}{dx} \right) &= 0, \\
\frac{2p}{w} \frac{dw}{dx} + \frac{1}{2} \frac{dp}{dx} + \frac{1}{2} \frac{p}{r} \frac{dr}{dx} &= 0, \qquad (\text{multiply with } w) \\
\frac{2pr}{w} \frac{dw}{dx} + \frac{1}{2} \left( r \frac{dp}{dx} + p \frac{dr}{dx} \right) &= 0,  \qquad (\text{divide with } pr)	\\
\frac{2}{w} \frac{dw}{dx} + \frac{1}{2} \frac{1}{pr} \frac{d}{dx} \left(p r \right) &= 0, \\
\frac{1}{w} \frac{dw}{dx} &= - \frac{1}{4} \frac{1}{pr} \frac{d}{dx} \left(p r \right), \\
\frac{d}{dx} \ln w &= \frac{d}{dx} \ln \left(pr\right)^{-1/4}, \\
w &= \left(pr\right)^{-1/4}.
\end{align*}
Before establishing the third relationship that involves the invariant function $I$, we need the following first-order and second-order derivative operators. The are similar to what we use at the beginning of the proof, albeit the role of $x$ and $t$ is reversed:
\begin{align*}
\frac{d}{dx} &= \frac{1}{\dot{x}} \frac{d}{dt}, \qquad \qquad \text{and} \\
\frac{d^2}{dx^2} &= \frac{d}{dx} \left(\frac{1}{\dot{x}} \right) \, \frac{d}{dt} + \frac{1}{\dot{x}} \frac{d}{dx} \left(\frac{d}{dt} \right) 
                  = \frac{1}{\dot{x}} \frac{d}{dt} \left(\frac{1}{\dot{x}} \right) \, \frac{d}{dt} + \frac{1}{\left(\dot{x} \right)^2} \, \frac{d^2}{dt^2} 
                  = \frac{1}{\left(\dot{x} \right)^2} \, \frac{d^2}{dt^2} - \frac{\ddot{x}}{\left(\dot{x} \right)^3} \frac{d}{dt}.
\end{align*}
Transforming the invariant function $I$ with an independent variable $x$ to the one that depends on $t$ gives the following:
\begin{align*}
I(x) &= \frac{q}{r} + \left[\frac{2}{w^2} \left(\frac{dw}{dx}\right)^2 - \frac{1}{w} \frac{d^2w}{dx^2} \right] \, \left(\frac{dx}{dt}\right)^2 - \frac{1}{w} \frac{dw}{dx} \frac{d^2x}{dt^2} \\
I(t) &= \frac{q}{r} + \left[\frac{1}{\left(\dot{x}\right)^2} \frac{2}{w^2} \left(\frac{dw}{dt}\right)^2 - \frac{1}{\left(\dot{x}\right)^2} \frac{1}{w} \frac{d^2w}{dt^2} + \frac{\ddot{x}}{\left(\dot{x}\right)^3} \frac{1}{w} \frac{dw}{dt} \right] \, \left(\frac{dx}{dt}\right)^2 - \frac{\ddot{x}}{\dot{x}} \frac{1}{w} \frac{dw}{dt} \\
     &= \frac{q}{r} + \frac{2}{w^2} \left(\frac{dw}{dt}\right)^2 - \frac{1}{w} \frac{d^2w}{dt^2} 	
      = \frac{q}{r} + w \left[ \frac{2}{w^3} \left(\frac{dw}{dt}\right)^2 - \frac{1}{w^2} \frac{d^2w}{dt^2} \right] \\
     &= \frac{q}{r} + w \frac{d}{dt} \left(- \frac{1}{w^2} \frac{dw}{dt} \right) 
      = \frac{q}{r} + w \frac{d}{dt} \left[\frac{d}{dt} \left(\frac{1}{w} \right) \right]
      = \frac{q}{r} + w \frac{d^2}{dt^2} \left(\frac{1}{w} \right).
\end{align*}
This completes the proof.
\end{proof}

\section{Reciprocal quadratic invariant function} 		\label{SectionPaine-n2}

In this section, we attempt to generalize the second Paine-de Hoog-Anderson (PdHA) problem, which is simply known in the literature as the second Paine problem~\cite{paine1981correction}. It is an SLP expressed in Schr\"odinger form, with Dirichlet boundary conditions. The corresponding generalized invariant function $I$ is given by a reciprocal binomial term with positive integer power, given as follows:
\begin{equation*}
I(t) = \frac{k}{\left(t + m\right)^n}, \qquad \text{where} \quad k, \; m > 0, \quad \text{and} \quad n \in \mathbb{N}.
\end{equation*}

In particular, we only focus on the case $n = 2$, but treating both the positive constants $k$ and $m$ as free parameters. A~discussion for other values of $n > 2$ will be presented in a separate work. This choice still yields the generalized second Paine problem with a reciprocal quadratic invariant function $I$, and leads to the following SLP in the Liouville normal form with Dirichlet boundary conditions:
\begin{equation}
-\frac{d^2v}{dt^2} + \frac{k}{\left(t + m \right)^2} v = \lambda \, v, \qquad v\left(0\right) = 0 = v\left(\pi\right).		\label{EquPaine-n2}
\end{equation}
In their original paper, Paine et al. took a special case of $k = 1$, $m = 0.1$, $n = 2$, $\alpha = 0$, and $\beta = \pi$~\cite{paine1981correction}. In what follows, we are interested in expressing SLP~\eqref{EquPaine-n2} in the canonical form.

\subsection{Vanishing potential and constant density functions}		\label{subsubsection-qvanish-rconstant}	

For this particular case, we have the following theorem.
\begin{theorem}		\label{Theorem-q0rconstant}
Let the potential function $q$ vanishes and the density function takes a constant value, that is, $q = 0$ and $r = r_0 \neq 0$, respectively. Then, for $k \neq 3/4$, the canonical form of ODE~\eqref{EquPaine-n2} is given as follows:
\begin{equation}
-\frac{d}{dx} \left\{ \accentset{\circ}{r}_0 \left[\left(2\rho + 1 \right) \left(x + x_0 \right)\right]^{4\rho/\left(2\rho + 1\right)} \frac{du}{dx} \right\} = \lambda \, r_0 \, u, \qquad a < x < b,	\label{EquPaine-rho0}
\end{equation}
where
\begin{align*}
\accentset{\circ}{r}_0 &= r_0^{\left(2\rho - 1\right)/\left(2\rho + 1\right)}, \\
2\rho + 1 &= 2 \pm \sqrt{1 + 4k }, \\
\frac{2\rho - 1}{2\rho + 1} &=  \frac{- \left(1 + 4k\right) \pm 2 \sqrt{1 + 4k}}{3 - 4k}, \\
\frac{4\rho}{2\rho + 1} &= \frac{2 \left(1 - 4k \pm \sqrt{1 + 4k}\right)}{3 - 4k}.
\end{align*}
For $k = 3/4$, ODE~\eqref{EquPaine-n2} may take one of the following two distinct canonical forms:
\begin{equation}
-\frac{d}{dx} \left\{ 8 \sqrt{r_0} \left(x + x_0 \right)^{3/2} \frac{du}{dx} \right\} = \lambda \, r_0 \, u, \qquad a < x < b,	\label{EquPaine-rho1}
\end{equation}
or
\begin{equation}
-\frac{d}{dx} \left\{ \frac{1}{r_0} e^{-2 r_0 \left(x + x_0 \right)} \frac{du}{dx} \right\} = \lambda \, r_0 \, u, \qquad a < x < b.		\label{EquPaine-rho2}
\end{equation}
All ODEs~\eqref{EquPaine-rho0}--\eqref{EquPaine-rho2} satisfy Dirichlet boundary conditions $u\left(a\right) = 0 = u\left(b\right)$. 
\end{theorem}

\begin{proof}
For the case of vanishing potential function $q = 0$ and nonzero constant density function $r = r_0 \neq 0$, we seek a function $w$ that satisfies the following ODE:
\begin{equation}
w \frac{d^2}{dt^2} \left(\frac{1}{w} \right) = \frac{k}{\left(t + m \right)^2}.		\label{EquPaine-n2a}
\end{equation}
Introducing a new dependent variable $\overline{\omega} = 1/w$, ODE~\eqref{EquPaine-n2a} can be written as follows, which turns out to be a special case of the Cauchy-Euler equation:
\begin{equation}
\left(t + m \right)^2 \frac{d^2\overline{\omega}}{dt^2} - k \overline{\omega} = 0.	\label{EquPaine-n2b}
\end{equation}
Introducing a new independent variable $\tau = t + m$, and because $d\tau = dt$ and $d^2\tau/dt^2 = 0$, we observe that the transformed ODE takes a similar form to ODE~\eqref{EquPaine-n2b}, that is, $\tau^2 \ddot{\overline{\omega}} - k \overline{\omega} = 0$, where double dots represent the second derivative with respect to $\tau$. Seeking an Ansatz in the form $\overline{\omega}\left(\tau \right) = \tau^{\rho}$, where $\tau = t + m$, we obtain $\rho^2 - \rho - k = 0$ as the indicial equation, which solves as
\begin{equation}
\rho = \rho_{1,2} = \frac{1}{2} \left(1 \pm \sqrt{1 + 4k} \right).		\label{EquPaine-indicial-roots}
\end{equation}
Observe that because $k > 0 > -1/4$, the indicial roots $\rho$ are always real-valued. The linearly independent solutions of ODE~\eqref{EquPaine-n2a} are thus given by
\begin{equation*}
w\left(\tau \right) = w_{1,2}\left(\tau \right) = \tau^{-\rho_{1,2}}.
\end{equation*}
From the Liouville transformation~\eqref{Equ-Liouville-transformation}, we can express the function $p$ and $dx/dt$ as follows:
\begin{equation}
p\left(\tau \right) = \frac{1}{r_0 \, w^4 \left(\tau \right)} = \frac{\overline{\omega}^4\left(\tau \right)}{r_0} = \frac{\tau^{4\rho}}{r_0},	\label{EquPaine-n2c}
\end{equation}
and 
\begin{equation}
\frac{dx}{d\tau} = \sqrt{\frac{p}{r_0}} = \frac{1}{r_0 \, w^2\left(\tau\right)} = \frac{\overline{\omega}^2\left(\tau\right)}{r_0} = \frac{\tau^{2\rho}}{r_0},		\label{EquPaine-n2d}
\end{equation}
because $dx/dt = dx/d\tau$. Integrating~\eqref{EquPaine-n2d} with respect to $\tau$, we obtain a relationship between the two independent variables $x$ and $t$:
\begin{equation}
x + x_0 = \frac{\tau^{2\rho + 1}}{r_0 \left(2\rho + 1 \right)}, \qquad x_0 \in \mathbb{R}.		\label{EquPaine-n2e}
\end{equation}
The Schr\"odinger variable $t$, or $\tau$, can be expressed explicitly in terms of the canonical variable $x$ using~\eqref{EquPaine-n2e}:
\begin{equation*}
t = - m + \left[r_0 \left(2\rho + 1\right) \left(x + x_0 \right)\right]^{1/\left(2\rho + 1 \right)}.
\end{equation*}
The function $p$ in~\eqref{EquPaine-n2c} can also be expressed in terms of the canonical variable $x$:
\begin{equation*}
p\left(x\right) = \frac{1}{r_0} \left[r_0 \left(2\rho + 1\right) \left(x + x_0 \right) \right]^{4\rho/\left(2\rho + 1\right)}.
\end{equation*}
Observe that the constant $r_0$ can be pulled out to form $\accentset{\circ}{r}_0$:
\begin{equation*}
\accentset{\circ}{r}_0 = r_0^{4\rho/\left(2\rho + 1\right) - 1} = r_0^{\left(2\rho - 1 \right)/\left(2\rho + 1\right)}.
\end{equation*}
Furthermore, using indicial roots $\rho$ in~\eqref{EquPaine-indicial-roots}, the three quantities involving $\rho$ stated in Theorem~\ref{Theorem-q0rconstant} can be verified using simple algebraic techniques of rationalizing the denominator:
\begin{align*}
2\rho + 1 &= 2 \pm \sqrt{1 + 4k }, \\
\frac{4\rho}{2\rho + 1} &= \frac{2 \pm 2\sqrt{1 + 4k}}{2 \pm \sqrt{1 + 4k}} = \frac{4 \pm 2 \sqrt{1 + 4k} - 2\left(1 + 4k \right)}{4 - \left(1 + 4k\right)} = \frac{2 \left(1 - 4k \pm \sqrt{1 + 4k}\right)}{3 - 4k}, \qquad k \neq \frac{3}{4}, \\
\frac{2\rho - 1}{2\rho + 1} &= \frac{\pm \sqrt{1 + 4k}}{2 \pm \sqrt{1 + 4k}} = \frac{\pm 2\sqrt{1 + 4k} - \left(1 + 4k \right)}{4 - \left(1 + 4k\right)} = \frac{-\left(1 + 4k \right) \pm 2 \sqrt{1 + 4k}}{3 - 4k}, \qquad k \neq \frac{3}{4}.
\end{align*}
To find the left- and right-endpoints of the variable $x$, we use~\eqref{EquPaine-n2e} that corresponds to $\alpha = 0$ and $\beta = \pi$, respectively. Thus,
\begin{equation*}
a = -x_0 + \frac{m^{2\rho + 1}}{r_0 \left(2\rho + 1\right)}, \qquad \qquad \text{and} \qquad \qquad
b = -x_0 + \frac{\left(\pi + m\right)^{2\rho + 1}}{r_0 \left(2\rho + 1\right)}.
\end{equation*}
By expressing the function $w$ (and $\overline{\omega}$) in terms of the canonical variable $x$, we can investigate the boundary conditions. Using~\eqref{EquPaine-n2e}, we have the following expression for $w$:
\begin{equation*}
w\left(x\right)	= \left[r_0 \left(2\rho + 1\right) \left(x + x_0 \right)\right]^{-\rho/\left(2\rho + 1 \right)}.
\end{equation*}
We also know that $\delta_1 = 0 = \gamma_1$, whereas $\delta_2 = 1 = \gamma_2$. Hence, we can calculate the values of $\delta_0$ and $\gamma_0$ by evaluating $w^2\left(x\right)$ at $x = a$ and $x = b$, respectively, which gives the following values
\begin{align*}
\delta_0 &= \frac{1}{w^2 \left(a\right)} = \overline{\omega}^2\left(a\right) = m^{2\rho} = m^{1 \pm \sqrt{1 + 4k}}, \qquad \text{and} \\
\gamma_0 &= \frac{1}{w^2 \left(b\right)} = \overline{\omega}^2\left(b\right) = \left(\pi + m \right)^{2\rho} = \left(\pi + m \right)^{1 \pm \sqrt{1 + 4k}}.
\end{align*}
However, because $\delta_1 = 0 = \gamma_1$, these quantities are irrelevant as we can always divide with each of them and the right-hand sides of the boundary conditions remain identical, which are Dirichlet-type. Thus, $u\left(a\right) = 0 = u\left(b\right)$. We have completed the proof for the first part of the theorem, that is, for the case $k \neq 3/4$. 

To prove the second part of the theorem, we take $k = 3/4$, and the indicial roots are given by
\begin{equation*}
\rho_1 = \frac{3}{2}, \qquad \qquad \text{and} \qquad \qquad \rho_2 = -\frac{1}{2}.
\end{equation*}
For the former, we have $2\rho_1 + 1 = 4$, $4\rho_1/\left(2\rho_1 + 1 \right) = 3/2$, and $\left(2 \rho_1 - 1\right)/\left(2\rho_1 + 1\right) = 1/2$. Consequently, the SLP~\eqref{EquPaine-n2} admits the following canonical form:
\begin{equation*}
-\frac{d}{dx} \left\{ 8 \sqrt{r_0} \left[\left(x + x_0 \right)\right]^{3/2} \frac{du}{dx} \right\} = \lambda \, r_0 \, u, \qquad a < x < b.		
\end{equation*}
The values of the boundary points $a$ and $b$ are given as follows:
\begin{equation*}
a = -x_0 + \frac{1}{4} \frac{m^4}{r_0}, \qquad \qquad \text{and} \qquad \qquad b = -x_0 + \frac{1}{4} \frac{\left(\pi + m \right)^4}{r_0}.
\end{equation*}
Irrelevant constant values $\delta_0 = m^3$ and $\gamma_0 = \left(\pi + m\right)^3$.

For the latter, we acquire
\begin{equation*}
\frac{dx}{d\tau} = \frac{1}{r_0 \, \tau},
\end{equation*}
and upon integration
\begin{equation*}
\tau = \pm e^{r_0 \left(x + x_0 \right)}, \qquad x_0 \in \mathbb{R}.
\end{equation*}
Moreover, because
\begin{equation*}
p\left(\tau \right) = \frac{1}{r_0 \, \tau^2},
\end{equation*}
we obtain
\begin{equation*}
p\left(x \right) = \frac{1}{r_0} e^{-2 r_0 \left(x + x_0\right)}.
\end{equation*}
Substituting this function $p$ to the SLP in the canonical form, we obtain the desired result, and the boundary points $a$ and $b$ are given as follows:
\begin{equation*}
a = -x_0 + \frac{1}{r_0} \ln m, \qquad \qquad \text{and} \qquad \qquad b = -x_0 + \frac{1}{r_0} \ln\left(\pi + m \right).
\end{equation*}
Irrelevant constant values $\delta_0 = 1/m$ and $\gamma_0 = 1/\left(\pi + m\right)$. Both ODEs~\eqref{EquPaine-rho1} and~\eqref{EquPaine-rho2} also admit Dirichlet boundary conditions $u\left(a\right) = 0 = u\left(b\right)$.
Thus, the proof is completed.
\end{proof}

\subsection{Constant potential and quadratic density functions}		

For the case of constant potential function $q = q_0 \neq 0$ and quadratic density function $r \left(\tau \right) = \tau^2$, we have three distinct subcases, depending on whether the indicial equation admits equal roots, distinct real roots, or complex conjugate roots. Each of these cases corresponds to the following relationships between $k$ and $q_0$, respectively:
\begin{itemize}[leftmargin=1em]
\item Case~A: $1 + 4 k = 4q_0$,
\item Case~B: $1 + 4 k > 4q_0$,
\item Case~C: $1 + 4 k < 4q_0$.
\end{itemize}

Following a similar derivation as in Subsubsection~\ref{subsubsection-qvanish-rconstant}, the function $\overline{\omega} = 1/w$ satisfies the following ODE:
\begin{equation}
\tau^2 \frac{d^2\overline{\omega}}{d\tau^2} + \left(q_0 - k \right) \overline{\omega} = 0,		\label{EquPaine-n2-rquadratic}
\end{equation}
with indicial equation $\rho^2 - \rho - \left(k - q_0 \right) = 0$ and its roots are given by
\begin{equation*}
\rho = \frac{1 \pm \sqrt{1 + 4 \left(k - q_0\right)}}{2}.
\end{equation*}
For Case~A, $\rho = 1/2$. The corresponding linear independent solutions to ODE~\eqref{EquPaine-n2-rquadratic} are given as follows:
\begin{equation*}
\overline{\omega}_1\left(\tau \right) = \sqrt{\tau}, \qquad \qquad \text{and} \qquad \qquad
\overline{\omega}_2\left(\tau \right) = \sqrt{\tau} \, \ln \tau, \qquad \tau > 0.
\end{equation*}

We have the following theorem.
\begin{theorem}[Equal roots~A1, part 1 of case A]
The SLP~\eqref{EquPaine-n2} admits the following canonical form with Dirichlet boundary conditions:
\begin{equation*}
-\frac{d^2u}{dx^2} + q_0 \, u = \lambda \, e^{2 \left(x + x_0\right)}\, u, \qquad a < x < b, \qquad u\left(a\right) = 0 = u\left(b\right),
\end{equation*}
where
\begin{align*}
a &= -x_0 + \ln m,  	\qquad 		  &b &= -x_0 + \ln \left(\pi + m\right), \\
\delta_0 &= m, 			\qquad &\gamma_0 &= \pi + m.
\end{align*}
\end{theorem}

\begin{proof}
For the case A1, using $\overline{\omega}_1\left(\tau\right) = \sqrt{\tau}$, $\tau > 0$, we observe that
\begin{equation*}
p\left(\tau \right) = \frac{\overline{\omega}_1^4}{r} = \frac{\tau^{2}}{\tau^2} = 1,
\end{equation*}	
and
\begin{equation*}
\frac{dx}{d\tau} = \sqrt{\frac{p}{r}} = \sqrt{\frac{1}{\tau^2}} = \frac{1}{\tau}.
\end{equation*}
Upon integration, we obtain
\begin{equation}
x + x_0 = \ln \tau, \qquad \qquad \text{or} \qquad \qquad \tau = e^{x + x_0}, \qquad x_0 \in \mathbb{R}.		\label{EquPaine-n2f}
\end{equation}
The functions $p$, $r$, and $w$ expressed in terms of the canonical variable $x$ are given as follows:
\begin{equation*}
p\left(x\right) = 1, \qquad \qquad r\left(x\right) = e^{2\left(x + x_0\right)}, \qquad \qquad \text{and} \qquad \qquad w\left(x\right) = e^{-\frac{1}{2}\left(x + x_0\right)}.
\end{equation*}
To find the values of $a$ and $b$, we substitute $t = \alpha = 0$ and $t = \beta = \pi$ to~\eqref{EquPaine-n2f}:
\begin{align*}
a + x_0 = \ln \left(0 + m\right)   \qquad &\Longrightarrow \qquad a = -x_0 + \ln m, \\
b + x_0 = \ln \left(\pi + m\right) \qquad &\Longrightarrow \qquad b = -x_0 + \ln \left(\pi + m\right).
\end{align*}
Finally, Dirichlet boundary conditions are confirmed if we can show that both $\delta_0$ and $\gamma_0$ are nonzero constants:
\begin{equation*}
\delta_0 = \overline{\omega}_1^2\left(a\right) = m \neq 0, \qquad \qquad \text{and} \qquad \qquad
\gamma_0 = \overline{\omega}_1^2\left(b\right) = \pi + m \neq 0.
\end{equation*}
The proof is completed.
\end{proof}

We have the following theorem for the case~A2.
\begin{theorem}[Equal roots~A2, part 2 of case A]
The SLP~\eqref{EquPaine-n2} admits the following canonical form: 
\begin{equation*}
-\frac{d}{dx} \left\{ \left[3 \left(x + x_0\right)\right]^{4/3} \frac{du}{dx} \right\} + q_0 \, u = \lambda \, e^{2 \left[3 \left(x + x_0 \right)\right]^{1/3}}  \, u, \qquad a < x < b,
\end{equation*}
with Dirichlet boundary conditions $u\left(a\right) = 0 = u\left(b\right)$, where
\begin{align*}
a &= -x_0 + \frac{1}{3} \ln^3 m,  	\qquad 		  &b &= -x_0 + \frac{1}{3} \ln^3 \left(\pi + m\right),\\
\delta_0 &= m \ln^2 m, 				\qquad &\gamma_0 &= \left(\pi + m \right) \ln \left(\pi + m\right).
\end{align*}
\end{theorem}

\begin{proof}
For the case A2, we utilize $\overline{\omega}_2 \left(\tau \right) = \sqrt{\tau} \ln \tau$, $\tau > 0$. We obtain
\begin{equation*}
p\left(\tau \right) = \frac{\overline{\omega}_2^4}{r} = \ln^4 \tau,
\end{equation*}
and
\begin{equation}
\frac{dx}{d\tau} = \sqrt{\frac{p}{r}} = \frac{\ln^2 \tau}{\tau}.		\label{EquPaine-n2g}
\end{equation}
We acquire an explicit expression for $x$ in terms of $\tau$ upon integrating~\eqref{EquPaine-n2g}:
\begin{equation}
x + x_0 = \int \frac{\ln^2 \tau}{\tau} \, d\tau = \frac{1}{3} \ln^3 \tau, \qquad x_0 \in \mathbb{R},	\label{EquPaine-n2h}
\end{equation}
or, an expression for $\tau$ in terms of $x$, that is,
\begin{equation*}
\tau = e^{\left[3 \left(x + x_0 \right)\right]^{1/3}}.
\end{equation*}
The functions $p$, $r$, and $w$ in terms of the canonical variable $x$ are given as follows:
\begin{align*}
p\left(x\right) &= \left[3 \left(x + x_0 \right)\right]^{\frac{4}{3}}, \\
r\left(x\right) &= e^{2\left[3 \left(x + x_0 \right)\right]^{\frac{1}{3}}},  \qquad \qquad \text{and} \\
w\left(x\right) &= \left[3 \left(x + x_0 \right)\right]^{-\frac{1}{3}} \, e^{-\frac{1}{2}\left[3 \left(x + x_0 \right)\right]^{\frac{1}{3}}}.
\end{align*}
The values of the endpoints $x = a$ and $x = b$ can be calculated using~\eqref{EquPaine-n2h} by substituting $t = \alpha = 0$ and $t = \beta = \pi$, respectively:
\begin{align*}
a + x_0 = \frac{1}{3} \ln^3 \left(0 + m \right) 	&\qquad \Longrightarrow \qquad a = -x_0 + \frac{1}{3} \ln^3 m, \\
b + x_0 = \frac{1}{3} \ln^3 \left(\pi + m\right) 	&\qquad \Longrightarrow \qquad b = -x_0 + \frac{1}{3} \ln^3 \left(\pi + m \right).
\end{align*}
Finally, Dirichlet boundary conditions are guaranteed when both coefficients $\delta_0$ and $\gamma_0$ are nonvanishing:
\begin{align*}
\delta_0 &= \overline{\omega}_2^2\left(a\right) = m \ln^2 m \neq 0, \\
\gamma_0 &= \overline{\omega}_2^2\left(b\right) = \left(\pi + m \right) \ln^2 \left(\pi + m \right) \neq 0.
\end{align*} 
The proof is completed.
\end{proof}

For Case~B, we have the following theorem.
\begin{theorem}[Case~B]
For $1 + 4k > 4q_0$, the SLP~\eqref{EquPaine-n2} satisfies the following canonical form: 
\begin{equation}
-\frac{d}{dx} \left[\left(2\rho - 1\right)^2 \left(x + x_0 \right)^2 \frac{du}{dx} \right] + q_0 \, u = \lambda \left[\left(2\rho - 1 \right) \left(x + x_0 \right) \right]^{2/\left(2\rho - 1\right)} \, u, \qquad a < x < b,	\label{Equ-CaseB}
\end{equation}
where
\begin{align*}
a = - x_0 + \frac{1}{2\rho - 1} m^{\left(2\rho - 1\right)}, \qquad &\qquad 
b = - x_0 + \frac{1}{2\rho - 1} \left(\pi + m\right)^{\left(2\rho - 1\right)}, \\
\delta_0 = m^{2\rho}, \qquad &\qquad \gamma_0 = \left(\pi + m\right)^{2 \rho}.
\end{align*}
The ODE~\eqref{Equ-CaseB} satisfies Dirichlet boundary conditions $u\left(a\right) = 0 = u\left(b\right)$. 
\end{theorem}

The proof for this theorem follows a similar argument as the proof of Theorem~\ref{Theorem-q0rconstant} presented in Subsubsection~\ref{subsubsection-qvanish-rconstant}. We only outline the main points.
\begin{proof}
Using $\overline{\omega} = \tau^{\rho}$, we have $p = \tau^{4\rho - 2}$ and $dx/d\tau = \tau^{\left(2\rho - 2\right)}$. After separating the variables and integrating with respect to each variable, we obtain
\begin{equation*}
x + x_0 = \frac{1}{2\rho - 1} \tau^{2\rho - 1}, \qquad \qquad \text{or} \qquad \qquad
\tau = \left[\left(2\rho - 1\right) \left(x + x_0 \right)\right]^{1/\left(2\rho - 1\right)}.
\end{equation*}
Observe that since $1 + 4k > 4q_0$, $\rho \neq 1/2$, and hence the denominator $2 \rho - 1$ never vanishes.
The functions $p$, $r$, and $w$ expressed in terms of the canonical variable $x$ are given by
\begin{align*}
p\left(x\right) &= \left(2\rho - 1\right)^2 \left(x + x_0 \right)^2, \\
r\left(x\right) &= \left[\left(2\rho - 1 \right) \left(x + x_0 \right) \right]^{2/\left(2\rho - 1\right)},  \\
w\left(x\right) &= \left[\left(2\rho - 1 \right) \left(x + x_0 \right) \right]^{-\rho/\left(2\rho - 1\right)}.
\end{align*}
The values of $a$ and $b$ can be calculated by substituting $\alpha = 0$ and $\beta = \pi$ to the relationship between $x$ and $t$, respectively. Similarly, the constants $\delta_0$ and $\gamma_0$ can be found by calculating $\overline{\omega}$ at $x = a$ and $x = b$, respectively. This completes the proof.
\end{proof}

For Case~C, the two roots of the indicial equation are complex conjugate, that is, $\rho = \frac{1}{2} \pm i \mu$, where $\mu = \frac{1}{2} \sqrt{4 \left(q_0 - k\right) - 1}$. The two linearly independent solutions to ODE~\eqref{EquPaine-n2-rquadratic} are given by
\begin{equation*}
\overline{\omega}_1 \left(\tau\right) = \sqrt{\tau} \, \cos \left(\mu \ln \tau \right), \qquad \text{and} \qquad
\overline{\omega}_2 \left(\tau\right) = \sqrt{\tau} \, \sin \left(\mu \ln \tau \right), \qquad \tau > 0.
\end{equation*}
We have the following theorem for the case corresponds to $\overline{\omega}_1$.
\begin{theorem}[Complex roots~C1, the first part of Case~C]
For $1 + 4k < 4q_0$, the SLP~\eqref{EquPaine-n2} could satisfy the following SLP in the canonical form in an asymptotic manner, with Dirichlet boundary conditions: 
\begin{equation*}
-\frac{d}{dx} \left\{\cos^4 \left[\mu \ln \left(1 + x + x_0 \right) \right] \frac{du}{dx} \right\} + q_0 \, u = \lambda \left(1 + x + x_0 \right) \, u, \qquad a < x < b, \qquad x_0 \in \mathbb{R},
\end{equation*}
where
\begin{align*}
a = -x_0 + m - 1, \qquad & \qquad 
b = -x_0 + \pi + m - 1, \\
\delta_0 = m \, \cos^2 \left(\mu \ln m \right), \qquad & \qquad 
\gamma_0 = \left(\pi + m\right) \, \cos^2 \left[ \mu \ln \left(\pi + m\right) \right], \qquad \text{and} \\
\mu \ln m \neq \pi \left(n - \frac{1}{2}\right), \qquad & \qquad \mu \ln \left(\pi + m\right) \neq \pi \left(n - \frac{1}{2}\right), \qquad n \in \mathbb{Z}.
\end{align*}
\end{theorem}

\begin{proof}
For case~C1, we acquire the following information
\begin{align}
p\left(\tau \right) &= \frac{\overline{\omega}_1^4}{r} = \cos^4 \left(\mu \ln \tau \right), \nonumber \\
\frac{dx}{d\tau} 	&= \sqrt{\frac{p}{r}} = \frac{\cos^2 \left(\mu \ln \tau \right)}{\tau}.		\label{EquPaine-n2i}
\end{align}	
Separating the variables in~\eqref{EquPaine-n2i} and integrating each side of the equation, we obtain
\begin{equation}
x + x_0 = \frac{1}{2} \ln \tau + \frac{1}{4 \mu} \sin \left(2 \mu \ln \tau \right), \qquad x_0 \in \mathbb{R}.				\label{EquPaine-n2j}
\end{equation}
By expanding the right-hand side of~\eqref{EquPaine-n2j}, let say about $\tau = 1$, we obtain
\begin{align*}
\ln \tau &= \left(\tau - 1 \right) - \frac{1}{2} \left(\tau - 1\right)^2 + \frac{1}{3} \left(\tau - 1\right)^3 + \cdots, \qquad 0 < \tau \leq 2, \\
\sin \left(2\mu \ln \tau \right) &= 2 \mu \left(\tau - 1 \right) - \mu \left(\tau - 1\right)^2 + \frac{2}{3} \mu \left(1 - 2\mu^2 \right) \left(\tau - 1\right)^3 + \cdots, \qquad \tau > 0, \\
x + x_0 &=  \left(\tau - 1\right) - \frac{1}{2}\left(\tau - 1\right)^2 + \frac{1}{3}\left(1 - \mu^2 \right) \left(\tau - 1\right)^3 + \cdots, \qquad 0 < \tau \leq 2.
\end{align*}
By considering the asymptotic expansion only up to the linear term, we can express $\tau$ explicitly in terms of the canonical variable $x$, given as follows (taking only the positive root):
\begin{equation*}
\tau = 1 + x + x_0.
\end{equation*}
The functions $p$, $r$, and $w$ can now be expressed in terms of $x$ asymptotically, given by
\begin{align*}
p\left(x\right) &= \cos^4 \left[\mu \ln \left(1 + x + x_0 \right) \right], \\
r\left(x\right) &= \left(1 + x + x_0 \right)^2, \\
w\left(x\right) &= \frac{\sec \left[\mu \ln \left(1 + x + x_0 \right) \right]}{\sqrt{1 + x + x_0}}.
\end{align*}
The values of $x = a$ and $x = b$ can be calculated using $x - x_1 = -\frac{1}{2}\left(t + m - 1 \right)^2$ by substituting $t = 0$ and $t = \pi$, respectively:
\begin{equation*}
a = -x_0 + m - 1, \qquad \qquad \text{and} \qquad \qquad b = -x_0 + \pi + m - 1.
\end{equation*}
Furthermore, by restricting $\mu \ln m \neq \pi \left(n - 1/2\right)$ and $\mu \ln \left(\pi + m\right) \neq \pi \left(n + 1/2\right)$, $n \in \mathbb{Z}$, it follows that both $\delta_0$ and $\gamma_0$ are nonzero, which guarantees that the associated ODE satisfies Dirichlet boundary conditions. These constants can be easily calculated and are given by
\begin{equation*}
\delta_0 = \overline{\omega}_1^2\left(a\right) = m \, \cos^2 \left(\mu \ln m \right), \qquad \qquad \text{and} \qquad \qquad 
\gamma_0 = \overline{\omega}_1^2\left(b\right) = \left(\pi + m\right) \, \cos^2 \left[ \mu \ln \left(\pi + m\right) \right].	
\end{equation*}
This completes the proof.
\end{proof}

We have the following theorem for the case corresponds to $\overline{\omega}_2$.
\begin{theorem}[Complex roots~C2, the second part of Case~C]
For $1 + 4k < 4q_0$, the SLP~\eqref{EquPaine-n2} may satisfy the following SLP in the canonical form in an asymptotic manner: 
\begin{equation*}
-\frac{d}{dx} \left( \sin^4 \left\{\mu \ln \left(1 + \sqrt[3]{\frac{3x}{\mu^2} \left[x + x_0 \right]} \, \right) \right\} \frac{du}{dx} \right) + q_0 \, u = \lambda \left(1 + \sqrt[3]{\frac{3x}{\mu^2} \left[x + x_0\right]} \, \right)^2 \, u, \qquad a < x < b,
\end{equation*}
where the ODE satisfies Dirichlet boundary conditions with
\begin{align*}
a = -x_0 + \frac{\mu^2}{3} \left(m - 1\right)^3, \qquad & \qquad
b = -x_0 + \frac{\mu^2}{3} \left(\pi + m - 1\right)^3, \qquad x_0 \in \mathbb{R},\\
\delta_0 = m \, \sin^2 \left(\mu \ln m \right), \qquad & \qquad 
\gamma_0 = \left(\pi + m \right) \, \sin^2 \left[ \mu \ln \left(\pi + m\right) \right], \qquad \text{and} \\
\mu \ln m \neq n \pi , \qquad & \qquad \mu \ln \left(\pi + m\right) \neq n \pi , \qquad n \in \mathbb{Z}.
\end{align*}
\end{theorem}

\begin{proof}
For case~C2, we obtain the following information
\begin{align}
p\left(\tau \right) &= \frac{\overline{\omega}_2^4}{r} = \sin^4 \left(\mu \ln \tau \right), \nonumber \\
\frac{dx}{d\tau} 	&= \sqrt{\frac{p}{r}} = \frac{\sin^2 \left(\mu \ln \tau \right)}{\tau}.		\label{EquPaine-n2k}
\end{align}	
Separating the variables in~\eqref{EquPaine-n2k}, integrating each side of the equation, we acquire
\begin{equation}
x + x_0 = \frac{1}{2} \ln \tau - \frac{1}{4 \mu} \sin \left(2 \mu \ln \tau \right), \qquad x_0 \in \mathbb{R}.	\label{EquPaine-n2l}
\end{equation}
Without loss of generality, we consider a Taylor-series expansion of the right-hand side of~\eqref{EquPaine-n2l} about $\tau = 1$ to obtain
\begin{align*}
x + x_0 &= \left[\frac{1}{2} \left(\tau - 1 \right) - \frac{1}{2} \left(\tau - 1\right)\right] + \left[-\frac{1}{4} \left(\tau - 1\right)^2 + \frac{1}{4} \left(\tau - 1\right)^2 \right] \\
		& \qquad + \left[\frac{1}{6} \left(\tau - 1\right)^3 - \frac{1}{6} \left(1 - 2\mu^2 \right) \left(\tau - 1\right)^3 \right] + {\cal O} \left(\tau - 1\right)^4, \\
		&= \frac{\mu^2}{3} \left(\tau - 1\right)^3 + {\cal O} \left(\tau - 1\right)^4. 
\end{align*}
By considering the asymptotic expansion only up to the cubic term, we can express $\tau$ explicitly in terms of~$x$, which is given by
\begin{equation*}
\tau = 1 + \sqrt[3]{\frac{3}{\mu^2} \left(x + x_0 \right)}.
\end{equation*}
The functions $p$, $r$, and $w$ expressed in terms of $x$ are given as follows:
\begin{align*}
p\left(x\right) &= \sin^4 \left\{\mu \ln \left(1 + \sqrt[3]{\frac{3}{\mu^2} \left[x + x_0\right]} \, \right) \right\}, \\
r\left(x\right) &= \left(1 + \sqrt[3]{\frac{3}{\mu^2} \left[x + x_0 \right]} \, \right)^2, \\
w\left(x\right) &= \frac{\csc \left\{\mu \ln \left(1 + \sqrt[3]{\frac{3}{\mu^2} \left[x + x_0\right]} \, \right) \right\}}{\sqrt{1 + \sqrt[3]{\frac{3}{\mu^2} \left[x + x_0\right]}}}.
\end{align*}
The values of $x = a$ and $x = b$ can be calculated using the relationship between the two independent variables by substituting $t = 0$ and $t = \pi$, respectively. It follows that
\begin{equation*}
a = - x_0 + \frac{\mu^2}{3} \left(m - 1\right)^3, \qquad \qquad \text{and} \qquad \qquad
b = - x_0 + \frac{\mu^2}{3} \left(\pi + m - 1\right)^3.
\end{equation*}
Finally, the constants $\delta_0$ and $\gamma_0$ can be calculated by substituting $x = a$ and $x = b$ to the squared of $\overline{\omega}_2$, respectively. They are given as follows:
\begin{align*}
\delta_0 = \overline{\omega}_2^2\left(a \right) = m \, \sin^2 \left(\mu \ln m \right), \qquad \qquad \text{and} \qquad \qquad
\gamma_0 = \overline{\omega}_2^2\left(b \right) = \left(\pi + m \right) \, \sin^2 \left[\mu \ln \left(\pi + m \right) \right].
\end{align*}
By imposing the restrictions $\mu \ln m \neq n \pi$ and $\mu \ln \left(\pi + m\right) \neq n \pi$, $n \in \mathbb{Z}$, these constants are nonzero, and thus the associated ODE admits Dirichlet boundary conditions.
This completes the proof.	
\end{proof}

\subsection{Both nonzero constant potential and density functions}		

We consider the case where both potential and density functions are nonzero constants, that is, $q = q_0 \neq 0$, and $r = r_0 \neq 0$. Consequently, the function $\overline{\omega} = 1/w$ satisfies a transformed version of the Bessel differential equation, where the first derivative of $\overline{\omega}$ with respect to the Schr\"odinger variable $\tau$ is absence.
\begin{equation}
\overline{\tau}^2 \frac{d^2\overline{\omega}}{d\overline{\tau}^2} + \left(\overline{\tau}^2 - k \right) \overline{\omega} = 0, \qquad \text{where} \quad \overline{\tau} = \sqrt{\frac{\left|q_0\right|}{\left|r_0\right|}} \tau.	\label{EquPaine-Bessel}
\end{equation}
Before we seek linearly independent solutions of ODE~\eqref{EquPaine-Bessel}, we state the following theorem on a transformation of the Bessel ODE, which we will need and use it later.
\begin{theorem}[Bowman's transformation~\cite{bowman1958introduction}]	\label{Theorem-Bowman}
Bowman (1958) gave a transformed version of the Bessel differential equation, given as follows:
\begin{equation}
\overline{\tau}^2 \frac{d^2\overline{\omega}}{d\overline{\tau}^2} + \left(2 \overline{p} + 1 \right) \frac{d\overline{\omega}}{d\overline{\tau}} + \left(\overline{\alpha}^2 \overline{\tau}^{2 \overline{r}} + \overline{\beta}^2 \right) \overline{\omega} = 0.	
\label{BowmanODE}	
\end{equation}
This transformed Bessel ODE~\eqref{BowmanODE} possesses the following solution
\begin{equation*}
\overline{\omega}\left(\overline{\tau }\right) = \frac{1}{\overline{\tau}^{\overline{p}}} \left[C_1 \, J_{\overline{q}/\overline{r}} \left(\frac{\overline{\alpha}}{\overline{r}} \overline{\tau}^{\overline{r}}\right) + 
C_2 \, Y_{\overline{q}/\overline{r}} \left(\frac{\overline{\alpha}}{\overline{r}} \overline{\tau}^{\overline{r}}\right)  \right],	\qquad C_1, \; C_2 \in \mathbb{R},
\end{equation*}
where $\overline{q} = \sqrt{\overline{p}^2 - \overline{\beta}^2}$, $J_n\left(\overline{\tau} \right)$ and $Y_n\left(\overline{\tau}\right)$, $n \in \mathbb{R}$, are the Bessel functions of the first and second kinds, respectively~\cite{gray1895a,watson1996a,korenev2002bessel,press2007numerical,nambudiripad2014bessel}.
\end{theorem}
Comparing ODEs~\eqref{EquPaine-Bessel} and~\eqref{BowmanODE}, we observe that $\overline{\alpha}^2 = 1$, $\overline{\beta}^2 = -k$, $\overline{p} = -1/2$, $\overline{q} = \frac{1}{2} \sqrt{4k + 1}$, and $\overline{r} = 1$. Hence, the associated linearly independent solutions of Bessel ODE~\eqref{EquPaine-Bessel} are given by
\begin{equation*}
\overline{\omega}_1\left(\overline{\tau} \right) = \sqrt{\overline{\tau}} \, J_{\frac{1}{2} \sqrt{4k + 1}} \left(\overline{\tau}\right), \qquad \text{and} \qquad
\overline{\omega}_2\left(\overline{\tau} \right) = \sqrt{\overline{\tau}} \, Y_{\frac{1}{2} \sqrt{4k + 1}} \left(\overline{\tau}\right).
\end{equation*}

We have the following theorem in connection to the Bessel function of the first kind $J_n$.
\begin{theorem}
The SLP~\eqref{EquPaine-n2} may satisfy the following canonical form asymptotically:
\begin{equation*}
-\frac{d}{dx} \left( \frac{4}{r_0} \, \overline{x}^2 \, J^4_{\frac{1}{2} \sqrt{4k + 1}} \left(2 \, \overline{x} \right) \frac{du}{dx} \right) + q_0 \, u = \lambda r_0 \, u, \qquad a < x < b,
\end{equation*}
where
\begin{align*}
\overline{x} &= \left[\frac{\Gamma_{\vartriangle}}{2} \sqrt{\left|q_0 r_0 \right|} \left(x + x_0\right) \right]^{1/\left(2 + \sqrt{4k + 1}\right)}, \qquad x_0 \in \mathbb{R},\\
\frac{1}{\Gamma_{\vartriangle}} &= \left(\frac{1}{\Gamma^2\left(1 + \frac{1}{2} \sqrt{4k + 1}\right)} - \frac{1}{\Gamma\left(\frac{1}{2} \sqrt{4k + 1}\right) \, \Gamma\left(2 + \frac{1}{2} \sqrt{4k + 1}\right)} \right), \\
a &= -x_0 + \frac{1}{\Gamma_{\vartriangle}} \, \frac{m}{r_0} \left(\frac{m}{2} \sqrt{\left|\frac{q_0}{r_0} \right|}\right)^{1 + \sqrt{4k + 1}}, \\
b &= -x_0 + \frac{1}{\Gamma_{\vartriangle}} \, \frac{\pi + m}{r_0} \left(\frac{\pi + m}{2} \sqrt{\left|\frac{q_0}{r_0} \right|}\right)^{1 + \sqrt{4k + 1}}, \\
\delta_0 &= \sqrt{\left|\frac{q_0}{r_0}\right|} m \; J^2_{\frac{1}{2}\sqrt{4k + 1}} \left(\sqrt{\left|\frac{q_0}{r_0}\right|} m \right), \\
\gamma_0 &= \sqrt{\left|\frac{q_0}{r_0}\right|} \left(\pi + m\right) \, J^2_{\frac{1}{2}\sqrt{4k + 1}} \left(\sqrt{\left|\frac{q_0}{r_0}\right|} \left(\pi + m\right) \right),
\end{align*}
$\Gamma$ denotes the gamma function, $J_n$ denotes the Bessel function of the first kind, and the terms $\sqrt{\left|q_0/r_0\right|} \, m$ and $\sqrt{\left|q_0/r_0\right|} \left(\pi + m\right)$ do not satisfy the zero $J_{\frac{1}{2}\sqrt{4k + 1}}$.
\end{theorem}

\begin{proof}
Let us consider the case of $\overline{\omega}_1$, where the function $p$ is given by the Bessel function of the first kind~$J_n$:
\begin{equation*}
p\left(\overline{\tau}\right) = \frac{\overline{\omega}_1^4}{r_0} = \frac{\overline{\tau}^2}{r_0} \, J_{\frac{1}{2} \sqrt{4k + 1}}^4 \left(\overline{\tau}\right).
\end{equation*}
The relationship between the canonical variable $x$ and Schr\"odinger variable $t$ is given by
\begin{equation*}
\sqrt{\frac{\left|q_0\right|}{\left|r_0\right|}} \, \frac{dx}{d\overline{\tau}} = \sqrt{\frac{p}{r}} = \frac{\overline{\tau}}{r_0} \, J_{\frac{1}{2} \sqrt{4k + 1}}^2 \left(\overline{\tau}\right).
\end{equation*}
It follows that
\begin{align}
x + x_0 = \int dx &= \sqrt{\frac{\left|r_0\right|}{\left|q_0\right|}} \int \frac{\overline{\tau}}{r_0} \, J_{\frac{1}{2} \sqrt{4k + 1}}^2 \left(\overline{\tau}\right) \, d\overline{\tau} \nonumber \\
&= \frac{\overline{\tau}^2}{2 r_0} \sqrt{\frac{\left|r_0\right|}{\left|q_0\right|}} \left[J_{\frac{1}{2} \sqrt{4k + 1}}^2 \left(\overline{\tau}\right) - 
J_{\frac{1}{2} \sqrt{4k + 1} - 1} \left(\overline{\tau}\right) \, J_{\frac{1}{2} \sqrt{4k + 1} + 1} \left(\overline{\tau}\right) \right].		\label{Equ-Paine-Bessel1}
\end{align}
For $0 \leq \overline{\tau} \ll 1$, the right-hand side of~\eqref{Equ-Paine-Bessel1} can be expressed asymptotically by
\begin{equation*}
x + x_0 = \frac{1}{\sqrt{\left|q_0 r_0 \right|}} \left(\frac{\overline{\tau}}{2}\right)^{\sqrt{4k + 1}} \left[ \frac{1}{\Gamma_{\vartriangle}} \, \frac{\overline{\tau}^2}{2} + {\cal O}\left(\overline{\tau}^4\right) \right],
\end{equation*}
where
\begin{equation*}
\frac{1}{\Gamma_{\vartriangle}} = \left(\frac{1}{\Gamma^2\left(1 + \frac{1}{2} \sqrt{4k + 1}\right)} - \frac{1}{\Gamma\left(\frac{1}{2} \sqrt{4k + 1}\right) \, \Gamma\left(2 + \frac{1}{2} \sqrt{4k + 1}\right)} \right).
\end{equation*}
Taking the lowest-order of asymptotic term, we now can express $\tau$ in terms of the canonical variable $x$:
\begin{equation*}
\tau = t + m = \Gamma_{\lozenge} \left(x + x_0 \right)^{1/\left(2 + \sqrt{4k + 1}\right)},
\end{equation*}
where
\begin{equation*}
\Gamma_{\lozenge} = 2 \sqrt{\left|\frac{r_0}{q_0}\right|} \left(\frac{\Gamma_{\vartriangle}}{2} \sqrt{\left|q_0 r_0\right|} \right)^{1/\left(2 + \sqrt{4k + 1}\right)}.
\end{equation*}
The functions $p$ and $w$ expressed in terms of the canonical variable $x$ are given asymptotically as follows, respectively:
\begin{align*}
p\left(x\right) &= \frac{4}{r_0} \, \overline{x}^2 \, J^4_{\frac{1}{2} \sqrt{4k + 1}} \left(2 \, \overline{x} \right), \\
w\left(x\right) &= \frac{1}{\sqrt{2 \, \overline{x}}} \, J^{-1}_{\frac{1}{2} \sqrt{4k + 1}} \left(2 \, \overline{x} \right), 
\end{align*}
where 
\begin{equation*}
\overline{x} = \left[\frac{\Gamma_{\vartriangle}}{2} \sqrt{\left|q_0 r_0 \right|} \left(x + x_0\right) \right]^{1/\left(2 + \sqrt{4k + 1}\right)}.	
\end{equation*}
When $t = 0$, the left canonical boundary is $x = a$, and both quantities satisfy the following relationship:
\begin{equation*}
m = \Gamma_{\lozenge} \left(a + x_0 \right)^{1/\left(2 + \sqrt{4k + 1}\right)} = 2 \sqrt{\left|\frac{r_0}{q_0}\right|} \left[\frac{\Gamma_{\vartriangle}}{2} \sqrt{\left|q_0 r_0\right|} \left(a + x_0 \right) \right]^{1/\left(2 + \sqrt{4k + 1}\right)},
\end{equation*}
for which the boundary $a$ can be expressed explicitly as follows:
\begin{align*}
\frac{\Gamma_{\vartriangle}}{2} \sqrt{\left|q_0 r_0 \right|} \left(a + x_0 \right)	&= \left(\frac{m}{2} \sqrt{\left|\frac{q_0}{r_0} \right|}\right)^{2 + \sqrt{4k + 1}}, \\
a &= -x_0 + \frac{1}{\Gamma_{\vartriangle}} \, \frac{m}{r_0} \left(\frac{m}{2} \sqrt{\left|\frac{q_0}{r_0} \right|}\right)^{1 + \sqrt{4k + 1}}.
\end{align*}
Similarly, when $t = \pi$, the right-hand side canonical boundary is $x = b$, and they satisfy the following relationship:
\begin{align*}
2 \sqrt{\left|\frac{r_0}{q_0}\right|} \left[\frac{\Gamma_{\vartriangle}}{2} \sqrt{\left|q_0 r_0\right|} \left(b + x_0 \right) \right]^{1/\left(2 + \sqrt{4k + 1}\right)} &= \pi + m, \\
\text{or} \qquad b &= -x_0 + \frac{1}{\Gamma_{\vartriangle}} \, \frac{\pi + m}{r_0} \left(\frac{\pi + m}{2} \sqrt{\left|\frac{q_0}{r_0} \right|}\right)^{1 + \sqrt{4k + 1}}.
\end{align*}
Finally, to ensure Dirichlet boundary conditions, we must ascertain that the terms $\sqrt{\left|q_0/r_0\right|} \, m$ and $\sqrt{\left|q_0/r_0\right|}$ $\left(\pi + m\right)$ do not satisfy as one of the zeros of $J_{\frac{1}{2}\sqrt{4k + 1}}$. This verifies that both $\delta_0$ and $\gamma_0$ are nonzero:
\begin{align*}
\delta_0 &= \overline{\omega}_1^2 \left(x = a \right) = \sqrt{\left|\frac{q_0}{r_0}\right|} m \; J^2_{\frac{1}{2}\sqrt{4k + 1}} \left(\sqrt{\left|\frac{q_0}{r_0}\right|} m \right) \neq 0, \\
\gamma_0 &= \overline{\omega}_1^2 \left(x = b \right) = \sqrt{\left|\frac{q_0}{r_0}\right|} \left(\pi + m\right) \, J^2_{\frac{1}{2}\sqrt{4k + 1}} \left(\sqrt{\left|\frac{q_0}{r_0}\right|} \left(\pi + m\right) \right) \neq 0.
\end{align*}
The proof is complete.
\end{proof}

We have the following theorem in connection to the Bessel function of the second kind $Y_n$.
\begin{theorem}
The SLP~\eqref{EquPaine-n2} could satisfy another canonical form asymptotically in connection to the Bessel function of the second kind $Y_n$, $n \geq 0$: 
\begin{equation*}
-\frac{d}{dx} \left\{ \pi^2 q_0 \left(x + x_1\right)^2 \, Y^4_{n/2} \left[\pi \sqrt{\left|q_0 r_0 \right|}  \left(x + x_1 \right) \right] \frac{du}{dx} \right\} + q_0 \, u = \lambda r_0 \, u, \qquad a < x < b,
\end{equation*}
where
\begin{align*}
a &= -x_1 + \frac{m}{\pi r_0}, \qquad x_1 \in \mathbb{R}, \qquad &\qquad 
b &= -x_1 + \frac{\pi + m}{\pi r_0}, \\	
\delta_0 &= \sqrt{\left|\frac{q_0}{r_0}\right|} m \, Y^2_{n/2} \left(\sqrt{\left|\frac{q_0}{r_0}\right|} m \right), \qquad & \qquad
\gamma_0 &= \sqrt{\left|\frac{q_0}{r_0}\right|} \left(\pi + m\right) \, Y^2_{n/2} \left(\sqrt{\left|\frac{q_0}{r_0}\right|} \left(\pi + m\right) \right),
\end{align*}
and the terms $\sqrt{\left|q_0/r_0\right|} \, m$ and $\sqrt{\left|q_0/r_0\right|} \left(\pi + m\right)$ must not satisfy the zero $Y_{n/2}$.
\end{theorem}

\begin{proof}
Similar as previously, we consider the case of $\overline{\omega}_2$, where now the function $p$ is expressible in terms of $Y_n$:
\begin{equation*}
p\left(\overline{\tau} \right) = \frac{\overline{\omega}_2^4}{r_0} = \frac{\overline{\tau}^2}{r_0} \, Y_{\frac{1}{2} \sqrt{4k + 1}}^4 \left(\overline{\tau}\right).
\end{equation*}
We have the following relationship between $x$ and $\overline{\tau}$
\begin{equation*}
\sqrt{\left|\frac{q_0}{r_0}\right|} \, \frac{dx}{d\overline{\tau}} = \sqrt{\frac{p}{r}} = \frac{\overline{\tau}}{r_0} \, Y_{\frac{1}{2} \sqrt{4k + 1}}^2 \left(\overline{\tau}\right).
\end{equation*}
After separating the variables and integrating both sides of the equation, it follows that
\begin{align}
x + x_1 = \int dx &= \sqrt{\left|\frac{r_0}{q_0}\right|} \int \frac{\overline{\tau}}{r_0} \, Y_{\frac{1}{2} \sqrt{4k + 1}}^2 \left(\overline{\tau}\right) \, d\overline{\tau} \nonumber \\
	&= \frac{\overline{\tau}^2}{2 r_0} \sqrt{\left|\frac{r_0}{q_0}\right|} \left[Y_{\frac{1}{2} \sqrt{4k + 1}}^2 \left(\overline{\tau}\right) - 
	Y_{\frac{1}{2} \sqrt{4k + 1} - 1} \left(\overline{\tau}\right) \, Y_{\frac{1}{2} \sqrt{4k + 1} + 1} \left(\overline{\tau}\right) \right].		\label{Equ-Paine-Bessel2}
\end{align}
A series expansion about $\overline{\tau} \gg 1$ for the right-hand side of~\eqref{Equ-Paine-Bessel2} exists when $\sqrt{4k + 1} \geq 0$, that is, 
\begin{equation*}
x + x_1 = \frac{1}{\pi r_0} \sqrt{\left|\frac{r_0}{q_0}\right|} \left[\overline{\tau} - \frac{k}{2 \overline{\tau}} + {\cal O} \left(\frac{1}{\overline{\tau}^3}\right) \right], \qquad	k \geq -\frac{1}{4}.	 
\end{equation*}
By considering only the linear term in $\overline{\tau}$, we can express $\tau$ in terms of $x$ asymptotically:
\begin{equation*}
\overline{\tau} = \pi \sqrt{\left|q_0 r_0 \right|} \left(x + x_1\right), \qquad \text{or} \qquad 
\tau = t + m = \pi r_0 \left(x + x_1\right), \qquad x_1 \in \mathbb{R}.
\end{equation*}
We can now express both $p$ and $w$ asymptotically as functions of the canonical variable $x$, which are given as follows, respectively:
\begin{align*}
p\left(x\right) &= \pi^2 q_0 \left(x + x_1\right)^2 \, Y^4_{\frac{n}{2}} \left(\pi \sqrt{\left|q_0 r_0 \right|}  \left(x + x_1 \right) \right), \qquad n \geq 0, \\ 
w\left(x\right) &= \frac{1}{\sqrt[4]{\left|q_0 r_0 \right|} \sqrt{\pi \left(x + x_1 \right)}} \, Y^{-1}_{\frac{n}{2}} \left(\pi \sqrt{\left|q_0 r_0 \right|}  \left(x + x_1 \right) \right), \qquad x + x_1 > 0. 
\end{align*}
The left and right canonical boundaries $x = a$ and $x = b$ can be calculated straightforwardly:
\begin{align*}
a = -x_1 + \frac{m}{\pi r_0}, \qquad \qquad \text{and} \qquad \qquad b = -x_1 + \frac{\pi + m}{\pi r_0}.
\end{align*}
To guarantee that the ODE admits Dirichlet boundary conditions, we need to ensure that the terms $\sqrt{\left|q_0/r_0\right|} \, m$ and $\sqrt{\left|q_0/r_0\right|} \left(\pi + m\right)$ must not satisfy as one of the zeros of $Y_n$, and thus verifying that both $\delta_0$ and $\gamma_0$ are nonzero:
\begin{align*}
\delta_0 &= \overline{\omega}_2^2 \left(x = a \right) = \sqrt{\left|\frac{q_0}{r_0}\right|} m \, Y^2_{\frac{n}{2}} \left(\sqrt{\left|\frac{q_0}{r_0}\right|} m \right) \neq 0, \\
\gamma_0 &= \overline{\omega}_2^2 \left(x = b \right) = \sqrt{\left|\frac{q_0}{r_0}\right|} \left(\pi + m\right) \, Y^2_{\frac{n}{2}} \left(\sqrt{\left|\frac{q_0}{r_0}\right|} \left(\pi + m\right) \right) \neq 0.
\end{align*}
The proof is completed.
\end{proof}

\subsection{Reciprocal linear function for $w\left(t\right)$}

In this subsection, we consider the case where the transformation function $w$ that appears in Lemma~\ref{lemmadependent}, Lemma~\ref{lemmacombined}, and equations~\eqref{Equ-invariant-function}--\eqref{Equ-Liouville-transformation} is a reciprocal linear function in the independent variable $t$. In other words, because $\overline{\omega} = 1/w$, this also means that the function $\overline{\omega}$ is assumed to be linear in $t$. As a consequence of this special case, the second term of the invariant function~\eqref{Equ-invariant-function} vanishes, and thus its first term, the quotient $q/r$, takes the form of the reciprocal binomial quadratic function in $t$. This special case of the second Paine problem was considered briefly by Ledoux and Ixaru et al.~\cite{ledoux2007study,ixaru1999slcpm}, but the derivation was notably absent in both works. As mentioned earlier in the introduction, we attempt to generalize this particular case and demonstrate the derivation in obtaining the $p$-, density, and potential functions. We also note that the title of this subsection could also be written as ``Reciprocal quadratic function for $w\left(x\right)$'' because the Schr\"odinger variable~$t$ is expressed as a quadratic function in the canonical variable $x$. 

We have the following theorem.
\begin{theorem}
For the particular case when the transformation function for the dependent variables takes the form of reciprocal linear function, the SLP~\eqref{EquPaine-n2} in the Liouville normal form can be transformed to the canonical form by preserving the type of the boundary conditions, i.e., Dirichlet. For constants $C_1 > 0$ and $x_0 \in \mathbb{R}$, the latter is given as follows:
\begin{equation*}
-\frac{d}{dx} \left[\frac{1}{8} C_1^3 \left(x + x_0\right)^3 \frac{du}{dx} \right] + \frac{k}{2} C_1^3 \left(x + x_0 \right) \, u = \lambda \left[\frac{1}{2} C_1 \left(x + x_0 \right) \right]^5 \, u, \qquad a < x < b, 
\end{equation*}
where
\begin{align*}
a &= -x_0 + 2 \sqrt{\frac{m}{C_1}}, \qquad &
b &= -x_0 + 2 \sqrt{\frac{\pi + m}{C_1}}, \\
\delta_0 &= \left(C_1 \, m \right)^2, \qquad &
\gamma_0 &= C_1^2 \left(\pi + m \right)^2.
\end{align*}
\end{theorem}

\begin{proof}
Let $w\left(t\right) = 1/\left(C_0 + C_1 t\right)$, where $C_0$ and $C_1$ are both positive constants, then it can be calculated straightforwardly that $w \frac{d^2}{dt^2} \left(1/w\right) = 0$. Thus, the invariant function reduces to
\begin{equation}
\frac{q}{r} = \frac{k}{\left(t + m\right)^2}.		\label{quotient-qr1}
\end{equation}
Furthermore, from 
\begin{equation}
pr = \frac{1}{w^4} = \left(C_0 + C_1 t\right)^4,		\label{product-pr}
\end{equation}
we encounter an underdetermined system, where we need to seek three unknown functions, that is $p$, $q$, and $r$, but we only possess two equations. Although we have another relationship $dx/dt = \sqrt{p/r}$, it does not really resolve the issue unless we impose further restriction to the relationship between the canonical variable $x$ and the Schr\"odinger variable $t$, which would help in making an educated guess for the three functions $p$, $q$, and $r$. Now let assume that $t$ is a function of $x$ quadratically, and let also the potential and density functions admit similar form, that is, linear in $t$, but with distinct rational powers, given as follows:
\begin{align*}
q\left(t\right) &= Q_0 \left(C_0 + C_1 t \right)^{1/2}, \qquad Q_0 > 0, \\
r\left(t\right) &= \left(C_0 + C_1 t \right)^{5/2}.
\end{align*}
The quotient $q/r$ becomes
\begin{equation}
\frac{q}{r} = \frac{Q_0/C_1^2}{\left(t + C_0/C_1\right)^2}. 		\label{quotient-qr2}	
\end{equation}
Comparing~\eqref{quotient-qr1} and~\eqref{quotient-qr2}, we obtain the following relationships:
\begin{equation*}
Q_0 = C_1^2 \, k, \qquad \qquad \text{and} \qquad \qquad C_0 = C_1 \, m.
\end{equation*}
Now, because $r$ is known, that is, $r\left(t\right) = \left[C_1 \left(t + m\right)\right]^{5/2}$, from the product~\eqref{product-pr}, we obtain an expression for $p$:
\begin{equation*}
p\left(t\right) = \left[C_1 \left(t + m \right)\right]^{3/2}.
\end{equation*}
We obtain the following relationship between the two independent variables:
\begin{equation}
\frac{dx}{dt} = \frac{1}{\sqrt{C_1 \left(t + m \right)}} \qquad \qquad \Longrightarrow \qquad \qquad x + x_0 = \int \frac{dt}{\sqrt{C_1 \left(t + m \right)}} = 2 \sqrt{\frac{t + m}{C_1}}, 	\label{equ-xt}
\end{equation}
where $x_0 \in \mathbb{R}$ is an integration constant. Because $1/w\left(t\right) = C_1 \left(t + m \right) = \left[C_1 \left(x + x_0 \right)/2 \right]^2 = 1/w\left(x\right)$, it follows that
\begin{align*}
p\left(x\right) &= \left[\frac{1}{2} C_1 \left(x + x_0 \right)\right]^3 = \frac{1}{8} C_1^3 \left(x + x_0 \right)^3, \\
q\left(x\right) &= C_1^2 k \, \left[\frac{1}{2} C_1 \left(x + x_0 \right)\right] = \frac{1}{2} k C_1^3 \left(x + x_0 \right), \\
r\left(x\right) &= \left[\frac{1}{2} C_1 \left(x + x_0 \right)\right]^5 = \frac{1}{32} C_1^5 \left(x + x_0 \right)^5.
\end{align*}
The canonical boundaries can be found by substituting $t = 0$ and $t = \pi$ to~\eqref{equ-xt}, respectively, and they are given as follows:
\begin{equation*}
a = -x_0 + 2 \sqrt{\frac{m}{C_1}}, \qquad \qquad \text{and} \qquad \qquad b = -x_0 + 2 \sqrt{\frac{\pi + m}{C_1}}.
\end{equation*}
Finally, Dirichlet boundary conditions can be confirmed by nonzero $\delta_0$ and $\gamma_0$:
\begin{align*}
\delta_0 &= \frac{1}{w^2\left(a\right)} = \left[\frac{1}{2} C_1 \left(a + x_0 \right)\right]^4 = C_1^2 \, m^2, \\
\gamma_0 &= \frac{1}{w^2\left(b\right)} = \left[\frac{1}{2} C_1 \left(b + x_0 \right)\right]^4 = C_1^2 \, \left(\pi + m \right)^2.
\end{align*}
This completes the proof.
\end{proof}

\begin{remark}
It is often practical to set the left-endpoint boundary $a = 0$, and thus, $x_0 = 2\sqrt{m/C_1}$, as it was done in~\cite{ledoux2007study,ixaru1999slcpm}. For this particular case, the right-endpoint boundary becomes
\begin{equation*}
b = \frac{2}{\sqrt{C_1}} \left(\sqrt{\pi + m} - \sqrt{m} \right) > 0.
\end{equation*}
The classical Paine problem takes $k = 1$ and $m = 0.1$, and by taking a special case $C_1 = 2$, we arrive to what Ledoux and Ixaru et al. stated~\cite{ledoux2007study,ixaru1999slcpm}, namely 
\begin{align*}
x_0 &= \sqrt{0.2}, \qquad \qquad &b &= \sqrt{2\pi + 0.2} - \sqrt{0.2}, \\
p\left(x\right) &= \left(x + \sqrt{0.2} \right)^3, \qquad &q\left(x\right) &= 4 \left(x + \sqrt{0.2}\right), &\qquad \text{and} \qquad r\left(x \right) = \left(x + \sqrt{0.2}\right)^5.
\end{align*}
\end{remark}

\begin{remark}
We can also generalize the powers of $q$ and $r$ to arbitrary positive numbers, let say, instead of $1/2$ and $5/2$, they become $n_q$ and $n_r$, respectively, where $n_r = n_q + 2$. It follows that 
\begin{equation*}
p\left(t\right) = \left[ C_1 \left(t + m \right) \right]^{4 - n_r},
\end{equation*}
and
\begin{equation*}
\frac{dx}{dt} = \left[ C_1 \left(t + m \right) \right]^{2 - n_r},
\end{equation*}
whereby upon integration, we obtain
\begin{equation*}
x + x_0 = \frac{C_1^{2 - n_r}}{3 - n_r} \left(t + m \right)^{3 - n_r}, \qquad n_r \neq 3.
\end{equation*}
Meanwhile, to avoid potential singularities for $p$, $q$, and $r$ in the canonical variable $x$, we impose a further restriction: because $n_q$ is positive, $n_r$ must be greater than 2. A similar argument applied to the relationship between $x$ and $t$ yields $n_r < 3$. This narrows the range for $n_r$ down to $2 < n_r < 3$. However, choosing $n_r$ values other than $5/2$ results in non-polynomial expressions for $p$, $q$, and $r$. While the case considered in this section is a special one, Ledoux and Ixaru et al.~\cite{ledoux2007study,ixaru1999slcpm} made a clever choice by taking $n_q = 1/2$ and $n_r = 5/2$, resulting in polynomial functions for $p$, $q$, and $r$ in the canonical variable.
\end{remark}

\section{Conclusion}		\label{SectionConclude}

We have considered the transformation of the Sturm-Liouville boundary value problem, often known as {SLP}, from its canonical form to the Schr\"odinger (Liouville normal) form and vice versa. Although it is theoretically possible to retrieve the SLP into its canonical form for any given SLP in the Schr\"odinger form, in practice, however, such an attempt is not always feasible. Implementing inverse Liouville's transformation can even be nearly impossible, in several, or even many, cases. 

For a particular case study, we investigated the second Paine-de Hoog-Anderson (PdHA) problem in a generalized manner. Also known in the literature as the Paine problem, the associated SLP with Dirichlet boundary conditions is given in its Liouville normal form instead of in its canonical form. The classical second Paine problem considered the corresponding invariant function in the form of a reciprocal binomial term with quadratic power and specific constants in the numerator and denominator, which are $1$ and $0.1$, respectively. We generalized these numbers to any positive constants while keeping the binomial power reciprocal quadratic.

Our study revealed that the difficulty of retrieving the SLP to the canonical form depends on the combinations between the potential and density functions. From the four special cases that we considered, inverting to SLP to the canonical form was relatively straightforward in some, while in others, such a process was impossible without adopting the technique of asymptotic expansion. One immediate consequence of the mentioned combinations is in the relationship between the independent variables, that is, whether we can easily find an exact expression for the Schr\"odinger variable in terms of the canonical variable by simply inverting the latter from the former. Whenever this fails, we simply proceed by finding its inverse asymptotically.

The exact SLP in its canonical form occurs in three cases, that is, when the potential function vanishes but the density function is nonzero constant, when the potential function is a nonzero constant and the density function is quadratic, and when the transformation function $w$ is reciprocal linear in $t$ (or, quadratic in $x$). For the second case, it occurs only in two special subcategories: when the roots of the indicial equations are either equal or real distinct. When these roots are complex conjugate, the SLP in the canonical form only accurate asymptotically. A similar case occurs when both the potential and density functions are nonzero constants. The $p$-functions appear in relatively elementary forms depending on the various cases and subcases. However, when both potential and density functions are nonzero constants, the $p$-function takes the form of the first and second kind Bessel functions. In all considered cases, Dirichlet boundary conditions follow accordingly.

A natural extension of this work would be to consider higher-order powers of the invariant function of the second Paine problem, such as reciprocal quartic, sextic, and even higher powers. This generalization, along with extending the reciprocal power to any positive real number, remains an open question. Further investigation could also explore other types of invariant functions found in various Liouville normal forms.

\subsection*{Conflicts of interest}
The authors declare no conflicts of interest. 

\subsection*{Acknowledgment}
The authors extend their sincere gratitude to Kenneth Liem from the Department of Chemical Engineering, Natural Science Campus, Sungkyunkwan University, for his stimulating and fruitful discussions.


\begin{thebibliography}{99}
\bibitem{zettl2010sturm} Zettl, A. (2010). \textit{Sturm-Liouville Theory}. Providence, Rhode Island, US: American Mathematical Society.	

\bibitem{luetzen1984sturm} Lützen, J. (1984). Sturm and Liouville's work on ordinary linear differential equations. The emergence of Sturm-Liouville theory. \textit{Archive for History of Exact Sciences} {\bfseries 29}: 309--376.

\bibitem{luetzen2012joseph} Lützen, J. (2012). Joseph Liouville 1809--1882: Master of Pure and Applied Mathematics. New York, US: Springer Science \& Business Media.
	
\bibitem{pryce1993numerical} Pryce, J. D. (1993). \textit{Numerical Solution of Sturm-Liouville Problems}. Oxford, UK: Oxford University Press.

\bibitem{bailey2001sleighn2} Bailey, P. B., Everitt, W. N., and Zettl, A. (2001). The SLEIGN2 Sturm-Liouville code. \textit{ACM Transaction on Mathematical Software} {\bfseries 27}(2): 143--192. 

\bibitem{amrein2005sturm} Amrein, W. O., Hinz, A. M., and Pearson, D. B. (Eds.). (2005). \textit{Sturm-Liouville Theory: Past and Present}. Basel, Switzerland: Birkhäuser Verlag.

\bibitem{gwaiz2008sturm} Al-Gwaiz, M. A. (2008). \textit{Sturm-Liouville Theory and Its Applications}. Berlin: Springer.	

\bibitem{agarwal2009ordinary} Agarwal, R. P., and O'Regan, D. (2009). \textit{Ordinary and Partial Differential Equations: With Special Functions, Fourier Series, and Boundary Value Problems}. Springer Science \& Business Media: New York, US.	

\bibitem{teschl2012ordinary} Teschl, G. (2012). \textit{Ordinary Differential Equations and Dynamical Systems}. Providence, Rhode Island, US: American Mathematical Society.

\bibitem{haberman2013applied} Haberman, R. (2013). \textit{Applied Partial Differential Equations with Fourier Series and Boundary Value Problems}, Fifth Edition. Boston, Massachusetts, US: Pearson Higher Education.	

\bibitem{guenther2019sturm} Guenther, R. B., and Lee, J. W. (2019). \textit{Sturm-Liouville Problems: Theory and Numerical Implementation}. Boca Raton, Florida, US: CRC Press.

\bibitem{kravchenko2020direct} Kravchenko, V. V. (2020). \textit{Direct and Inverse Sturm-Liouville Problems: A Method of Solution}. Cham, Switzerland: Birkhäuser Verlag, Springer Nature.

\bibitem{masjed-jamei2020special} Masjed-Jamei, M. (2020). \textit{Special Functions and Generalized Sturm-Liouville Problems}. Basel, Switzerland: Birkhäuser.

\bibitem{zettl2021recent} Zettl, A. (2021). \textit{Recent Developments in Sturm-Liouville Theory}. Berlin, Germany and Boston, Massachusetts, US: De Gruyter.

\bibitem{almdallal2009efficient} Al-Mdallal, Q. M. An efficient method for solving fractional Sturm--Liouville problems. \textit{Chaos Solitons Fractals} \textbf{2009}, \emph{40}, 183--189.

\bibitem{klimet2013fractional} Klimek, M. and Agrawal, O. P. Fractional Sturm--Liouville problem. \textit{Computers \& Mathematics with Applications} \textbf{2013}, \emph{66}, 795--812.

\bibitem{zayernouri2013fractional} Zayernouri, M. and Karniadakis, G. E. Fractional Sturm--Liouville eigen-problems: Theory and numerical approximation. \textit{Journal of Computational Physics} \textbf{2013}, \emph{252}, 495--517.

\bibitem{chanane2007computing} Chanane, B. (2007). Computing the spectrum of non-self-adjoint Sturm--Liouville problems with parameter-dependent boundary conditions. \textit{Journal of Computational and Applied Mathematics} {\bfseries 206}(1): 229--237.

\bibitem{veliev2007non} Veliev, O. A. (2007). Non-self-adjoint Sturm-Liouville operators with matrix potentials. \textit{Mathematical Notes} {\bfseries 81}(3-4): 440--448.

\bibitem{albeverio2008on} Albeverio, S., Hryniv, R., and Mykytyuk, Y. (2008). On spectra of non-self-adjoint Sturm--Liouville operators. \textit{Selecta Mathematica} {\bfseries 13}: 571--599.

\bibitem{behrndt2008accumulation} Behrndt, J., Katatbeh, Q., and Trunk, C. (2008). Accumulation of complex eigenvalues of indefinite Sturm--Liouville operators. \textit{Journal of Physics A: Mathematical and Theoretical} {\bfseries 41}(24): 244003.

\bibitem{xie2013non} Xie, B., and Qi, J. (2013). Non-real eigenvalues of indefinite Sturm--Liouville problems. \textit{Journal of Differential Equations} {\bfseries 255}(8): 2291--2301.

\bibitem{levitin2015accumulation} Levitin, M., and Seri, M. (2015). Accumulation of complex eigenvalues of an indefinite Sturm--Liouville operator with a shifted Coulomb potential. arXiv preprint arXiv:1503.08615.

\bibitem{sager1984the} Sager, H. C. (1984). The Sturm-Liouville equation with time-dependent boundary conditions. \textit{Journal of Mathematical Analysis and Applications} {\bfseries 102}(1): 275--287.

\bibitem{binding1994sturm} Binding, P. A., Browne, P. J., and Seddighi, K. (1994). Sturm--Liouville problems with eigenparameter dependent boundary conditions. \textit{Proceedings of the Edinburgh Mathematical Society} {\bfseries 37}(1): 57--72.

\bibitem{stikonas2007the} Štikonas, A. (2007). The Sturm-Liouville problem with a nonlocal boundary condition. \textit{Lithuanian Mathematical Journal} {\bfseries 47}: 336--351.

\bibitem{adam2017rays} Adam, J. (2017). \textit{Rays, Waves, and Scattering: Topics in Classical Mathematical Physics}. Princeton, New Jersey, US: Princeton University Press.

\bibitem{birkhoff1989ordinary} Birkhoff, G., and Rota, G.-C. (1989). \textit{Ordinary Differential Equations}, Fourth Edition. New York, US: John Wiley \& Sons.

\bibitem{everitt1982on} Everitt, W. N. (1982). On the transformation theory of ordinary second-order linear symmetric differential expressions. \textit{Czechoslovak Mathematical Journal} {\bfseries 32}(2): 275--306.

\bibitem{paine1981correction} Paine, J. W., de Hoog, F. R., and Anderssen, R. S. (1981). On the correction of finite difference eigenvalue approximations for Sturm-Liouville problems. \textit{Computing} {\bfseries 26}(2): 123--139.

\bibitem{ledoux2007study} Ledoux, V. (2007). \textit{Study of Special Algorithms for Solving Sturm-Liouville and Schr\"odinger Equations}. Doctoral dissertation. Ghent, Belgium: Ghent University.

\bibitem{ixaru1999slcpm} Ixaru, L. G., De Meyer, H., and Berghe, G. V. (1999). SLCPM12--A program for solving regular Sturm—Liouville problems. \textit{Computer Physics Communications} {\bfseries 118}(2-3): 259--277.

\bibitem{everitt2005a} Everitt, W. N. (2005). A catalogue of Sturm-Liouville differential equations. In Amrein, W. O., Hinz, A. M., and Pearson, D. B. (Eds.). \textit{Sturm-Liouville Theory: Past and Present}, pp.~271--331. Basel, Switzerland: Birkhäuser Verlag.

\bibitem{karjanto2022perturbed} Karjanto, N. (2022). Perturbed potential temperature field in the atmospheric boundary layer. \textit{ZAMM‐Journal of Applied Mathematics and Mechanics/Zeitschrift für Angewandte Mathematik und Mechanik} {\bfseries 102}(8): e202100484.

\bibitem{karjanto2022modified} Karjanto, N. (2022). On modified second Paine--de Hoog--Anderssen boundary value problem. \textit{Symmetry} {\bfseries 14}(1): 54.

\bibitem{arnold2019computing} Arnold, D. N., David, G., Filoche, M., Jerison, D., and Mayboroda, S. (2019). Computing spectra without solving eigenvalue problems. \textit{SIAM Journal on Scientific Computing} {\bfseries 41}(1): B69--B92.

\bibitem{filoche2016universal} Filoche, M. and Mayboroda, S. (2012). Universal mechanism for Anderson and weak localization. \textit{Proceedings of the National Academy of Sciences} {\bfseries 109}(37): 14761--14766.

\bibitem{liouville1837second} Liouville, J. (1837). Second mémoire sur le développement des fonctions ou parties de fonctions en séries dont les divers termes sont assujettis à satisfaire à une même équation différentielle du second ordre, contenant un paramètre variable. (Second dissertation on the development of functions or parts of functions in series whose various terms are subject to satisfying the same second order differential equation, containing a variable parameter.) \textit{Journal de Mathématiques Pures et Appliquées (Journal of Pure and Applied Mathematics)} {\bfseries 2}: 16--35. (in French)

\bibitem{bowman1958introduction} Bowman, F. (1958). \textit{Introduction to Bessel Functions}. New York, US: Dover.

\bibitem{gray1895a} Gray, A., and Mathews, G. B. (1895). \textit{A Treatise on Bessel Functions and Their Applications to Physics}. London, England, UK: Macmillan and Co. 

\bibitem{watson1996a} Watson, G. N. (1996). \textit{A Treatise on the Theory of Bessel Functions}, Second Edition. Cambridge, England, UK: Cambridge University Press.

\bibitem{korenev2002bessel} Korenev, B. G. (2002). \textit{Bessel Functions and Their Applications}. Boca Raton, Florida, US: CRC Press.

\bibitem{press2007numerical} Press, W. H., Teukolsky, S. A., Vetterling, W. T., \& Flannery, B. P. (2007). \textit{Numerical Recipes: The Art of Scientific Computing}, Third Edition. New York, US: Cambridge University Press.

\bibitem{nambudiripad2014bessel} Nambudiripad, K. B. M. (2014). \textit{Bessel Functions}. Oxford, England, UK: Alpha Science International Limited.
\end{thebibliography}
\end{document}